\pdfoutput=1

\documentclass{amsart}
\usepackage[numbers]{natbib}
\usepackage{mathrsfs}

\usepackage{bbm} 

\usepackage{microtype} 
\usepackage{braket}
\usepackage{adjustbox}

\usepackage{tikz-cd}
\tikzset{
    labl/.style={anchor=south, rotate=90, inner sep=.5mm}
}
\usetikzlibrary{decorations.markings}
\usetikzlibrary{intersections}
\usepackage{pgf,tikz}
\usetikzlibrary{arrows}

\definecolor{qqqqff}{rgb}{0,0,1}

\makeatletter
\tikzcdset{
  open/.code     = {\tikzcdset{hook, circled};},
  closed/.code   = {\tikzcdset{hook, slashed};},
  open'/.code    = {\tikzcdset{hook', circled};},
  closed'/.code  = {\tikzcdset{hook', slashed};},
  circled/.code  = {\tikzcdset{markwith = {\draw (0,0) circle (.375ex);}};},
  slashed/.code  = {\tikzcdset{markwith = {\draw[-] (-.4ex,-.4ex) -- (.4ex,.4ex);}};},
  markwith/.code ={
    \pgfutil@ifundefined%
    {tikz@library@decorations.markings@loaded}%
    {\pgfutil@packageerror{tikz-cd}{You need to say %
      \string\usetikzlibrary{decorations.markings} to use arrows with markings}{}}{}%
    \pgfkeysalso{/tikz/postaction = {
      /tikz/decorate,
      /tikz/decoration={markings, mark = at position 0.5 with {#1}}}
    }
  },
}
\makeatother
\usepackage[english]{babel}
\usepackage{csquotes}
\usepackage{verbatim}

\usepackage[all]{xy}

\usepackage{multicol}

\usepackage{lmodern}
\usepackage[T1]{fontenc}

\usepackage{amssymb}
\usepackage{amsmath}
\usepackage{amsthm}
\usepackage{mathtools}

\usepackage{todonotes} 

\usepackage{enumerate}

\usepackage{csquotes}
\usepackage{verbatim}

\usepackage[all]{xy}

\xyoption{2cell}
\UseAllTwocells

\usepackage{multicol}

\usepackage{lmodern}
\usepackage[T1]{fontenc}

\usepackage{xr-hyper}
\usepackage{hyperref}
\hypersetup{
colorlinks=true,
linkcolor=[RGB]{0,0,204},
filecolor=magenta,      
urlcolor=cyan,
citecolor=[RGB]{0,204,0},
}

\usepackage[nameinlink]{cleveref}
\usepackage[normalem]{ulem}
\newtheorem{theorem}{Theorem}[section]
\newtheorem{proposition}[theorem]{Proposition}
\newtheorem{corollary}[theorem]{Corollary}
\newtheorem{lemma}[theorem]{Lemma}

\theoremstyle{definition}
\newtheorem{definition}[theorem]{Definition}
\newtheorem{example}[theorem]{Example}

\newtheorem{claim}[theorem]{Claim}
\newtheorem{remark}[theorem]{Remark}

\numberwithin{equation}{subsection}

\def\lim{\mathop{\mathrm{lim}}\nolimits}

\def\Spec{\mathop{\mathrm{Spec}}}
\def\Hom{\mathop{\mathrm{Hom}}\nolimits}

\def\Pic{\mathop{\mathscr{P}\mathrm{ic}}\nolimits}
\def\pic{\mathop{\mathrm{Pic}}\nolimits}

\def\deg{\mathop{\mathrm{deg}}\nolimits}

\def\min{\mathop{\mathrm{min}}\nolimits}

\def\dim{\mathop{\mathrm{dim}}\nolimits}
\def\rk{\mathop{\mathrm{rk}}\nolimits}

\def\lim{\mathop{\mathrm{lim}}\nolimits}

\def\Aut{\mathop{\mathrm{Aut}}\nolimits}

\newcommand{\GL}{\mathrm{GL}} 
\newcommand{\cO}{\mathcal{O}}

\newcommand{\cV}{\mathcal{V}}
\newcommand{\cL}{\mathcal{L}}

\newcommand{\SCoh}{\mathscr{C}oh}
\newcommand{\cE}{\mathcal{E}}
\newcommand{\cW}{\mathcal{W}}
\newcommand{\cF}{\mathcal{F}}

\usepackage{graphicx} 

\def\Bun{\mathop{\mathscr{B}un}\nolimits}

\def\Spec{\mathop{\mathrm{Spec}}\nolimits}
\def\Bun{\mathop{\mathscr{B}un}\nolimits}

\def\lim{\mathop{\mathrm{lim}}\nolimits}

\def\Spec{\mathop{\mathrm{Spec}}}
\def\Hom{\mathop{\mathrm{Hom}}\nolimits}

\def\deg{\mathop{\mathrm{deg}}\nolimits}

\def\min{\mathop{\mathrm{min}}\nolimits}

\def\dim{\mathop{\mathrm{dim}}\nolimits}
\def\rk{\mathop{\mathrm{rk}}\nolimits}

\def\lim{\mathop{\mathrm{lim}}\nolimits}

\newcommand{\isom}[1][]{
	\ar[#1]
	\ar@<0.7ex>@{}[#1]|-*=0[@]{\sim}}

\usepackage{pbox}

\title{There are no exotic compact \\ moduli of sheaves on a curve}
\author{Andres Fernandez Herrero, Dario Wei{\ss}mann, and Xucheng Zhang}

\address{(A. Fernandez Herrero) University of Pennsylvania, Department of Mathematics,
209 South 33rd Street,
Philadelphia, PA 19104, USA,
Email: \href{mailto: andresfh@sas.upenn.edu}{andresfh@sas.upenn.edu}}
\address{(D. Wei{\ss}mann) Instytut Matematyczny Polskiej Akademii Nauk, ul. Śniadeckich 8, 00-656 Warszawa, Poland,
Email: \href{mailto: dario.weissmann@posteo.de}{dario.weissmann@posteo.de}}
\address{(X. Zhang) Tsinghua University, Yau Mathematical Sciences Center, Beijing 100084, China, 
Email: \href{mailto: withboundary@163.com}{withboundary@163.com}}

\date{}

\begin{document}

\begin{abstract}
We study moduli of coherent sheaves of some given degree
and positive rank on a curve. 
We show that there is only one nonempty open condition 
on families of sheaves that yields a universally closed adequate moduli space, 
namely, the one that recovers the classical moduli of slope semistable vector bundles.
\end{abstract}
\maketitle

\tableofcontents

\section{Introduction}
Let $C$ be a smooth projective connected curve over an algebraically closed
field $k$ of arbitrary characteristic. Assume that the genus of $C$ is at least $2$. 
In this paper we investigate the moduli problem of coherent sheaves on $C$.
There are two discrete numerical invariants attached to a coherent sheaf on $C$, 
namely, its rank and degree. 
Given a fixed pair of integers $(d,r)$,
coherent sheaves on $C$ of rank $r$ and degree $d$ are parametrized by a
connected algebraic stack $\SCoh_r^d$ which is locally of finite type over $k$. 

The case when the rank $r=0$ corresponds to studying torsion sheaves of a given length $d$ on $C$. This situation is simple from the point of view of moduli theory; in this case it is known that the stack $\SCoh_0^d$ is of finite type and admits a projective adequate moduli space, which coincides with the Hilbert scheme $\mathrm{Hilb}^d(C)$ of length $d$ subchemes of $C$ (we recall this in \Cref{prop: moduli of torsion sheaves}).

The case when the rank $r>0$ is more complex. In this case, the stack
$\SCoh^d_r$ is not of finite type, and it does not admit an adequate moduli
space. The classical way to deal with this problem is to restrict our attention
to a special class of sheaves which satisfy a stability condition. 
In order to truly capture a large collection of sheaves and 
have control over the deformation theory, 
such a stability condition is usually assumed to be
open in families. A more intrinsic way of thinking about such a
stability condition is to fix the choice of an open substack $\mathscr{U}
\subset \SCoh_r^d$. 

The most classical open substack of sheaves, defined by Mumford
\cite{mumford-projetive-invariants} and Seshadri \cite{seshadri-unitary} is the
substack $\Bun_r^{d,ss} \subset \SCoh_r^d$ of slope semistable vector bundles of
rank $r$ and degree $d$. It is known that $\Bun_r^{d,ss}$ admits an adequate
moduli space $M_r^d$, which is a projective normal variety over $k$. The notion
of slope semistability is natural from the point of view of differential
geometry \cite{narasimhan-seshadri-unitary, donaldson-narasimhan-seshadri} and
Geometric Invariant Theory (GIT) \cite{mumford-projetive-invariants,
mumford_git}. However, it is reasonable to wonder whether there are other exotic
stability conditions on $\SCoh^d_r$ which
yield proper moduli spaces. In this paper we answer this in the negative.
\begin{theorem}[\Cref{prop: moduli spaces in coh(C)} + \Cref{thm: proper moduli space is semistable locus}] \label{thm: main thm 1 intro}
    Let $C$ be a smooth projective connected curve of genus at least $2$ over an algebraically closed field $k$. Fix a pair of integers $(d,r)$ such that $r>0$. Then, the stack of semistable vector bundles $\Bun_r^{d,ss} \subset \SCoh_r^d$ is the only nonempty open substack which admits a universally closed adequate moduli space in the sense of \cite[Def. 5.1.1]{ams}.
\end{theorem}

We note that in \Cref{thm: main thm 1 intro} the moduli space is allowed to be a
highly non-separated algebraic space a priori. 
So we are asserting that the only possible
nonempty moduli space of positive rank sheaves on a curve satisfying a very mild
``compactness'' condition is the classical moduli of slope semistable vector
bundles.
\begin{remark}
We observe that \Cref{thm: main thm 1 intro} actually holds for geometrically connected smooth projective curves over an arbitrary
ground field (cf. \citep[Lem. 4.8]{weissmann-zhang}), because a vector bundle is
semistable if and only if its base-change to an algebraic closure is semistable
(see \citep[Cor. 1.3.8]{huybrechts-lehn}).
\end{remark}
\begin{remark}
    The hypothesis that $\mathscr{U}$ is universally closed is quite important; otherwise it is possible to obtain exotic moduli spaces of sheaves which are not schemes. Indeed, \citep[\S 5]{weissmann-zhang} exhibits an open substack $\mathscr{U} \subset \Bun^2_3$ containing unstable vector bundles that admits a non-schematic separated good moduli space. In this paper we provide an example of an open substack $\mathscr{U} \subset \Bun_3^{0,ss}$ which also admits a non-schematic separated good moduli space. This type of behavior cannot happen if $r=1$ or $r=2$ (see \Cref{thm: moduli spaces of rank 2 vector bundles}).
\end{remark}
As a direct consequence of \Cref{thm: main thm 1 intro}, we obtain the following result, which answers a question posed in \citep[Rmk. 5.12]{weissmann-zhang}.
\begin{corollary}
Let $\mathscr{U} \subset \SCoh_r^d$ be an open substack that admits an adequate moduli space. Then there exists an open substack $\overline{\mathscr{U}} \subset \SCoh_r^d$ containing $\mathscr{U}$ that admits a proper adequate moduli space if and only if $\mathscr{U} \subset \Bun_r^{d,ss}$.\qed
\end{corollary}

The techniques we employ in this paper are fundamentally different from the ones
in previous works by the authors \cite{weissmann-zhang,
distinguishing_spaces_paper}. In those papers, the main tools took advantage of
the strict restrictions imposed on an open substack that has a separated or
schematic moduli space. Indeed, the results \citep[Thm. 2.1]{weissmann-zhang}
and \citep[Thm. 5.11]{distinguishing_spaces_paper} show that if an open substack
$\mathscr{U} \subset \Bun_r^d$ has a separated schematic adequate moduli space
$U$, then $\mathscr{U}$ is contained in the semistable locus $\Bun_r^{d,ss}$ and
the induced morphism $U \to M_r^d$ is an open immersion. Furthermore, when
$r=2$, it is shown in \cite[Thm. 4.1]{weissmann-zhang} that $\mathscr{U} \subset
\Bun_2^{d,ss}$ under the weaker hypothesis that $U$ is a separated algebraic
space. The fact that the induced morphism $U \to M^d_2$ is an open immersion
remains true, and is proven in \Cref{thm: moduli spaces of rank 2 vector
bundles}.

Instead, in this paper we focus on compactness (i.e. universal closedness), and we allow the larger flexibility of non-separated and non-schematic moduli spaces. The tools we employ in this paper are more classical in spirit: the main ingredients in our proofs involve Luna's stratification of the moduli of semistable vector bundles, Harder-Narasimhan stratifications of the stack of vector bundles, and the resolution of Lange's conjecture. 

In addition to these classical tools, we make use of some results on sheaves with geometrically reductive automorphism groups which we derive in \Cref{section: sheaves with reductive automorphism groups}. An example is the following characterization of sheaves with geometrically reductive automorphism groups on a proper scheme, which could be of independent interest.
\begin{theorem}[ = \Cref{thm: main}]\label{thm: main thm 2 in intro}
    Let $X$ be a proper scheme over an algebraically closed field $k$. Let $\mathcal{E}$ be a coherent sheaf on $X$ such that its group of automorphisms $\Aut(\mathcal{E})$ is geometrically reductive. Then, we have $\mathcal{E} \cong \oplus_{i=1}^{\ell} \mathcal{E}_i^{\oplus n_i}$, where $\mathcal{E}_i$ is simple and $\Hom(\mathcal{E}_i, \mathcal{E}_j) =0$ for $i \neq j$.
\end{theorem}

\Cref{thm: main thm 2 in intro} is related to the geometry of the stack
$\SCoh(X)$ of coherent sheaves on $X$, in two aspects. First, in characteristic
0 the local structure theorems for algebraic stacks \cite[Thm.
1.2]{luna-slice-for-stacks} tell us that \'{e}tale locally around any closed
point with a reductive stabilizer the stack is a quotient stack of the form
$[\mathrm{Spec}(A)/G]$, where $G$ is a reductive group. Second, \cite[Prop.
9.3.4]{ams} implies that any closed point of a stack with an adequate moduli
space has a geometrically reductive stabilizer.

\medskip

\noindent \textbf{Acknowledgements.} The question of finding open substacks of bundles with proper moduli spaces was proposed by Jochen Heinloth to Xucheng Zhang during his Ph.D. studies. We would like to thank Jochen Heinloth for the constant and generous support on this topic. We would also like to thank Georg Hein, Andr\'es Ib\'a\~nez N\'u\~nez, Svetlana Makarova, and Jianping Wang for helpful mathematical discussions.
Dario Weißmann was partially funded by PRIN 2022
``Geometry of Algebraic Structures: Moduli, Invariants, Deformations''.

\section*{Notation and conventions} 
In this paper, we work over a fixed algebraically closed ground field $k$.
Unless otherwise stated, we denote by $C$ a 
smooth projective connected curve over $k$ of genus at least $2$. 
For a nonzero vector bundle $\mathcal{E}$ on $C$, 
its slope $\mu(\mathcal{E})$ is defined to be the quotient 
$\deg(\mathcal{E})/\rk(\mathcal{E})$ of its degree by its rank. 

A nonzero vector bundle $\mathcal{E}$ on $C$ is called (slope) semistable 
(resp. stable) if for all subbundles 
$0 \subsetneq \mathcal{F} \subsetneq \mathcal{E}$ 
we have  $\mu(\mathcal{E}) \geq \mu(\mathcal{F})$ 
(resp. $\mu(\mathcal{E}) > \mu(\mathcal{F})$). 
We say that a vector bundle is polystable if it is a direct sum of stable vector bundles
of the same slope.

We denote by $\Bun_r^d$ the stack of vector bundles of rank $r$ and degree $d$ on $C$. We denote by $\Bun_r^{d,ss} \subset \Bun_r^d$ the open substack of semistable vector bundles. We write $M_r^d$ for the adequate moduli space of $\Bun_r^{d,ss}$. There is an open substack $\Bun_r^{d,s} \subset \Bun_r^{d,ss}$ parameterizing stable vector bundles; we denote by $M_r^{d,s}$ the good moduli space of $\Bun_r^{d,s}$. 

The moduli space $M_r^d$ is a normal projective irreducible variety over $k$, which contains the smooth moduli space of stable vector bundles $M_r^{d,s}\subset M_r^d$ as a dense open subscheme. This is well-known from the constructions of these moduli spaces using GIT, see for example \citep[Prem. partie, III, Thm. 17]{seshadri}.

\section{Sheaves with reductive automorphism groups} \label{section: sheaves with reductive automorphism groups}
 Let $X$ be a proper scheme over $k$. In this section, we study coherent sheaves on $X$ with geometrically reductive automorphism groups.

\subsection{Polysimple sheaves}

By \cite[Thm. 2]{MR0086358}, the category of coherent sheaves on $X$ is a Krull-Schmidt category, i.e., every coherent sheaf $\mathcal{E}$ on $X$ is isomorphic to a direct sum
\[
\mathcal{E} \cong \mathcal{E}_1^{\oplus n_1} \oplus \cdots \oplus \mathcal{E}_{\ell}^{\oplus n_{\ell}},
\]
where the $\mathcal{E}_i$ are distinct indecomposable sheaves. The collection of pairs $(\mathcal{E}_i,n_i)$ is uniquely determined up to permutation. We will be concerned with sheaves that have a special type of Krull-Schmidt decomposition. We call such sheaves polysimple, in analogy to polystable vector bundles, which are direct sums of stable vector bundles. For the following definition, recall that a sheaf $\mathcal{E}$ is called simple if $\mathrm{End}(\mathcal{E})=k$ (as $k$-vector spaces).
\begin{definition}[Polysimple sheaf] \label{defn: polysimple sheaf}
A coherent sheaf $\mathcal{E}$ on $X$ is \emph{polysimple} if it is a direct sum $\mathcal{E} \cong \oplus_{i=1}^{\ell} \mathcal{E}_i^{\oplus n_i}$ of simple sheaves $\mathcal{E}_i$ and $\mathrm{Hom}(\mathcal{E}_i,\mathcal{E}_j)=0$ for $i \neq j$.
\end{definition}
If $\mathcal{E}$ is a polysimple sheaf with the direct sum decomposition as in \Cref{defn: polysimple sheaf}, then we have an isomorphism of algebraic groups
\[
\mathrm{Aut}(\mathcal{E}) \cong \mathrm{GL}_{n_1} \times \cdots \times \mathrm{GL}_{n_\ell},
\]
so in particular $\Aut(\mathcal{E})$ is (geometrically) reductive. Our goal in this section is to show that any coherent sheaf on $X$ with a geometrically reductive automorphism group is necessarily of this form.

Let us first collect a few standard results on indecomposable coherent sheaves, see \cite[Thm. 2, Lem. 6 and 7]{MR0086358}. For any coherent sheaf $\mathcal{E}$ on $X$, we denote by $\mathrm{End^{nil}}(\mathcal{E}) \subset \mathrm{End}(\mathcal{E})$ the subset of nilpotent endomorphisms of $\mathcal{E}$. 
\begin{lemma}\label{lem: end-of-indecomp}
Let $X$ be a proper scheme over $k$, and let $\mathcal{E}$ be an indecomposable coherent sheaf on $X$. Then the subset $\mathrm{End^{nil}}(\mathcal{E}) \subset \mathrm{End}(\mathcal{E})$ is a nilpotent two-sided ideal and 
\[
\mathrm{End}(\mathcal{E})/\mathrm{End^{nil}}(\mathcal{E})=k.
\]
In particular, we have $\Aut(\mathcal{E})$ is geometrically reductive if and only if $\mathrm{End^{nil}}(\mathcal{E})=0$, if and only if $\mathcal{E}$ is simple.
\end{lemma}
\begin{proof}
We will repeatedly use the fact that an endomorphism of an indecomposable sheaf is either nilpotent or is an automorphism, see \cite[Lem. 6]{MR0086358}. 

For any $a,b \in \mathrm{End^{nil}}(\mathcal{E})$, the endomorphism $a-b$ is not an automorphism as otherwise at least one of $a,b$ is an automorphism by \cite[Lem. 7]{MR0086358}, so we have $a-b \in \mathrm{End^{nil}}(\mathcal{E})$. 
Further, for any $a \in \mathrm{End^{nil}}(\mathcal{E})$ and 
$b \in \mathrm{End}(\mathcal{E})$, 
the compositions $a \circ b, b \circ a$ are not
automorphisms, so $a \circ b, b \circ a \in \mathrm{End^{nil}}(\mathcal{E})$.
This shows that $\mathrm{End^{nil}}(\mathcal{E}) \subset
\mathrm{End}(\mathcal{E})$ is a two-sided ideal and hence equals to the
(Jacobson) radical of $\mathrm{End}(\mathcal{E})$ (see, e.g. \cite[Chap. XVII,
Thm. 6.1 (b)]{Lang-alg}). To show its nilpotency, note that
$\mathrm{End}(\mathcal{E})=\Gamma(X,\mathcal{E} \otimes \mathcal{E}^*)$ is a
finite-dimensional $k$-algebra, its radical is always nilpotent (see, e.g.
\cite[Chap. XVII, Thm. 6.1 (d)]{Lang-alg}). Finally, as
$\mathrm{End}(\mathcal{E})/\mathrm{End^{nil}}(\mathcal{E})$ is a
finite-dimensional division algebra over the algebraically closed field $k$, it
follows that $\mathrm{End}(\mathcal{E})/\mathrm{End^{nil}}(\mathcal{E})=k$.
\end{proof}
\begin{theorem}\label{thm: main}
Let $X$ be a proper scheme over $k$. Then a coherent sheaf $\mathcal{E}$ on $X$ has a geometrically reductive automorphism group if and only if it is polysimple.
\end{theorem}
\begin{proof}
The if part is clear. For the only if part, let $\mathcal{E} \cong \oplus_{i=1}^\ell
\mathcal{E}_i^{\oplus n_i}$ be the decomposition of $\mathcal{E}$ into distinct
indecomposables. Set $n:=n_1+\cdots+n_\ell$. Each endomorphism of $\mathcal{E}$ can
be naturally viewed as a $n \times n$ matrix with entries in
$\Hom(\mathcal{E}_i,\mathcal{E}_j)$ and the composition of endomorphisms is
compatible with matrix multiplication. Let
\[
\mathrm{U}:=\bigoplus_{i \neq j} \mathrm{Hom}(\mathcal{E}_i^{\oplus n_i},\mathcal{E}_j^{\oplus n_j}) \oplus \bigoplus_{i=1}^\ell \mathrm{End^{nil}}(\mathcal{E}_i)^{\oplus n_i} \subset \mathrm{End}(\mathcal{E}).
\]
By \cite[Lem. 2.4]{indecom-Higgs}, this subset $\mathrm{U}$ is exactly the (Jacobson) radical of $\mathrm{End}(\mathcal{E})$. In particular, it is a two-sided nilpotent ideal in $\mathrm{End}(\mathcal{E})$. The subset $\mathrm{Id} + \mathrm{U} \subset \mathrm{Aut}(\mathcal{E})(k)$ 
is the set of $k$-points of an algebraic subgroup $W \subset \mathrm{Aut}(\mathcal{E})$. Indeed, for all $a,b \in \mathrm{U}$ with $b^m=0$ we have 
\begin{gather*}
(\mathrm{Id}+a) \circ (\mathrm{Id}+b)^{-1}=(\mathrm{Id}+a) \circ \sum_{i=0}^{m-1} (-1)^ib^i=\mathrm{Id}+\left[a+\sum_{i=1}^{m-1} (-1)^i(\mathrm{Id}+a) \circ b^i\right].
\end{gather*}
It is evident that $W$ is a smooth unipotent algebraic group, and it is normal because we have $g \circ a \circ g^{-1} \in \mathrm{U}$ for all $g \in \mathrm{Aut}(\mathcal{E})(k)$ and $a \in \mathrm{U}$. If the vector space $\mathrm{U}$ was non-zero, then $W$ would be a positive dimensional unipotent normal subgroup of $\mathrm{Aut}(\mathcal{E})$, contradicting the assumption that $\mathrm{Aut}(\mathcal{E})$ is geometrically reductive. Therefore, $\mathrm{Hom}(\mathcal{E}_i,\mathcal{E}_j)=0$ for $i \neq j$ and $\mathrm{End^{nil}}(\mathcal{E}_i)=0$ for all $i$. Since $\mathcal{E}_i$ is indecomposable, by \Cref{lem: end-of-indecomp} this implies that $\mathcal{E}_i$ is simple.
\end{proof}

We end this subsection with a lemma that we employ in later sections.
\begin{lemma}\label{thm:non-split}
Let $X$ be a scheme over $k$. For any non-split short exact sequence of coherent sheaves on $X$
\[
0 \to \mathcal{E}_1 \xrightarrow{i} \mathcal{E} \xrightarrow{p} \mathcal{E}_2 \to 0
\]
such that $\mathcal{E}_1,\mathcal{E}_2$ are simple and $\mathrm{Hom}(\mathcal{E}_1,\mathcal{E}_2)=0$, we have $\mathcal{E}$ is simple.
\end{lemma}
\begin{proof}

For any $\phi \in \mathrm{End}(\mathcal{E})$, the composition $p \circ \phi \circ i: \mathcal{E}_1 \to \mathcal{E}_2$ is zero by assumption, so $\phi \circ i$ factors through some morphism $\mathcal{E}_1 \to \mathcal{E}_1$, which is $\lambda \cdot \mathrm{Id}_{\mathcal{E}_1}$ for some scalar $\lambda \in k$ since $\mathcal{E}_1$ is simple. 
Consider the morphism $\phi-\lambda \cdot \mathrm{Id}_{\mathcal{E}}:
\mathcal{E} \to \mathcal{E}$. We conclude by showing that it is zero. Note that it is zero on $\mathcal{E}_1$ and thus factors through some
morphism $h: \mathcal{E}_2 \to \mathcal{E}$. It remains to show $h=0$ since $p: \mathcal{E} \twoheadrightarrow \mathcal{E}_2$ is surjective. 
To see this we consider the composition 
$p \circ h: \mathcal{E}_2 \to \mathcal{E} \to \mathcal{E}_2$. 
If it was non-zero, then it would be an isomorphism since
$\mathcal{E}_2$ is simple and $\mathcal{E}$ would be the split extension, a contradiction.
\end{proof}
\subsection{The case of a curve}
In this subsection we collect some results on coherent sheaves on curves which we use in later sections. We use $C$ to denote a smooth projective connected curve over $k$.

\begin{corollary} \label{lemma: geometrically reductive automorphism S-equivalent to direct sum of same stable bundle}
Suppose that a semistable vector bundle $\mathcal{E}$ on $C$ is S-equivalent to $\mathcal{E}_0^{\oplus n}$ for some stable vector bundle $\mathcal{E}_0$ and integer $n \geq 1$. Then $\mathrm{Aut}(\mathcal{E})$ is geometrically reductive if and only if $\mathcal{E} \cong \mathcal{E}_0^{\oplus n}$.
\end{corollary}
\begin{proof}
The if part is clear. Suppose $\mathrm{Aut}(\mathcal{E})$ is geometrically
reductive. Then by \Cref{thm: main} we have $\mathcal{E} \cong \oplus_{i=1}^\ell
\mathcal{E}_i^{\oplus n_i}$ for some simple bundles $\mathcal{E}_i$. Since
$\mathcal{E}$ is semistable and S-equivalent to $\mathcal{E}_0^{\oplus n}$, each
direct summand $\mathcal{E}_i$ is semistable and S-equivalent to
$\mathcal{E}_0^{\oplus m_i}$ for some integer $m_i \geq 1$. If $m_i \geq 2$ for
some $i \in \{1,\ldots,\ell\}$, then $\mathcal{E}_i$ admits a subbundle $i:
\mathcal{E}_0 \hookrightarrow \mathcal{E}$ and a quotient $p: \mathcal{E}_i
\twoheadrightarrow \mathcal{E}_0$. However, the non-zero morphism $i \circ p \in
\mathrm{End}(\mathcal{E}_i)$ cannot be a multiple of
$\mathrm{Id}_{\mathcal{E}_i}$, a contradiction.
\end{proof}
Recall that a torsion sheaf on $C$ is a coherent sheaf whose stalk at the generic point of $C$ is zero, therefore it is of rank $0$.
\begin{lemma}\label{cor: simple torsion sheaf on curve}
Any simple torsion sheaf on $C$ is of the form $\mathcal{O}_C/\mathfrak{m}_x$, where $\mathfrak{m}_x$ is the ideal sheaf of some closed point $x \in C(k)$.
\end{lemma}
\begin{proof}
Note that any simple torsion sheaf $\mathcal{E}$ on $C$ is supported on a single
closed point, say $x \in C(k)$, and hence is of the form
$\mathcal{O}_C/\mathfrak{m}_x^{n}$ for some integer $n \geq 1$. 
If $n \geq 2$, then $\mathcal{E}$ admits a subsheaf 
$i: \mathcal{O}_C/\mathfrak{m}_x \hookrightarrow \mathcal{E}$ 
and a quotient $p: \mathcal{E} \twoheadrightarrow \mathcal{O}_C/\mathfrak{m}_x$. 
The non-zero morphism $i \circ p \in \mathrm{End}(\mathcal{E})$ 
cannot be a multiple of $\mathrm{Id}_{\mathcal{E}}$, a contradiction. 
\end{proof}
 We have an explicit description of coherent sheaves on $C$ whose automorphism groups are geometrically reductive. 
\begin{proposition}\label{prop: torsion is non-reductive}
Let $\mathcal{E}$ be a coherent sheaf on $C$. Then $\mathrm{Aut}(\mathcal{E})$ is geometrically reductive if and only if one of the following holds:
\begin{enumerate}
\item 
$\mathcal{E}$ is a torsion sheaf of the form
\[
\mathcal{E} \cong \oplus_{i=1}^{\ell} \left(\mathcal{O}_{C}/\mathfrak{m}_{x_i}\right)^{\oplus n_i}
\]
for some closed points $x_1,\ldots,x_\ell \in C(k)$ and positive integers $n_i>0$.
\item 
$\mathcal{E}$ is a vector bundle of the form
\[
\mathcal{E} \cong \oplus_{i=1}^{\ell} \mathcal{E}_i^{\oplus n_i}
\]
for some simple vector bundles $\mathcal{E}_i$ and $\mathrm{Hom}(\mathcal{E}_i,\mathcal{E}_j)=0$ for $i \neq j$.
\end{enumerate}
In particular, a coherent sheaf on $C$ with torsion cannot have a geometrically reductive automorphism group unless it is a torsion sheaf.
\end{proposition}
\begin{proof}
The if part is clear. Suppose $\mathrm{Aut}(\mathcal{E})$ is geometrically
reductive. Then by \Cref{thm: main} we have $\mathcal{E} \cong \oplus_{i=1}^\ell
\mathcal{E}_i^{\oplus n_i}$ for some simple sheaves $\mathcal{E}_i$ with
$\mathrm{Hom}(\mathcal{E}_i,\mathcal{E}_j)=0$ for $i \neq j$. Note that any
coherent sheaf on a curve (non-canonically) decomposes into a direct sum of its
torsion-free quotient and torsion subsheaf. Thus, each $\mathcal{E}_i$ is either
a torsion sheaf or a vector bundle. Since there always exist non-zero morphisms
from vector bundles to torsion sheaves, the condition
$\mathrm{Hom}(\mathcal{E}_i,\mathcal{E}_j)=0$ for $i \neq j$ implies that
$\mathcal{E}_i$ are either all torsion sheaves or all vector bundles. In the
former case, we are done by \Cref{cor: simple torsion sheaf on curve}. In the
latter case, we obtain (2). 
\end{proof}

\section{Universally closed moduli of sheaves on a curve}

The main goal of this section is to complete the proof \Cref{thm: main thm 1 intro}. The proof proceeds via several steps. 

\begin{itemize}
    \item In \Cref{section:reduction to vector bundles} we reduce from an open substack of coherent sheaves
to an open substack of vector bundles via the restriction imposed on closed points
having to have geometrically reductive automorphism groups. This step is quite general and does not yet require the universally closed assumption or the assumption that the curve has genus at least $2$.

\item Secondly, in \Cref{section:reduction to semistable vector bundles} we reduce
to an open substack of the semistable locus. 
By considering a minimal Harder-Narasimhan stratum we obtain another obstruction to having geometrically reductive automorphism groups at closed points. 
This uses in an essential way the universally closed assumption as well as
that semistable vector bundles of different slopes admit non-zero morphisms between them (see \Cref{lemma:general-bundles-have-homs}), which is an immediate consequence of Lange's conjecture proven
by Russo and Teixidor i Bigas \cite{langesconjecture-1}.

\item Lastly, we exclude the possibility of a nontrivial open substack of the
semistable locus admitting a universally closed adequate moduli space via its
interaction with the Luna stratification of $M^d_r$. 
We recall this stratification in 
\Cref{section: lunas stratitifcation for the moduli of bundles} and
conclude the proof in \Cref{section:substacks of the semistable locus}.
\end{itemize}

\subsection{Reduction to bundles}
\label{section:reduction to vector bundles}
For this section we keep our choice of smooth projective connected curve $C$ over $k$.
Denote by $\mathscr{C}oh$ the stack of coherent sheaves on $C$. 
It is an algebraic stack with affine diagonal and locally of finite type over $k$ 
(see \cite[\href{https://stacks.math.columbia.edu/tag/09DS}{Tag 09DS},
\href{https://stacks.math.columbia.edu/tag/0DLY}{Tag 0DLY},
\href{https://stacks.math.columbia.edu/tag/0DLZ}{Tag 0DLZ}]{sp}). 
In this section we study open substacks of $\mathscr{C}oh$ that admit adequate
moduli spaces in the sense of \cite[Def. 5.1.1]{ams}.

Given a nonnegative integer $r \geq 0$ and an integer $d$, there is a connected
open and closed substack $\mathscr{C}oh_r^d \subset \mathscr{C}oh$ parametrizing
coherent sheaves on $C$ of rank $r$ and degree $d$. There is an open
substack $\mathscr{B}un_r^d \subset \mathscr{C}oh_r^d$ parametrizing locally
free sheaves on $C$, which is non-empty if $r>0$.

\begin{proposition}\label{prop: moduli spaces in coh(C)}
If $r>0$, then any quasi-compact open substack of $\mathscr{C}oh_r^d$ that admits an adequate moduli space is contained in $\mathscr{B}un_r^d$.
\end{proposition}
\begin{proof}
Let $\mathscr{U} \subset \mathscr{C}oh_r^d$ be a quasi-compact open substack
that admits an adequate moduli space. If $\mathcal{E} \in |\mathscr{U}|$ has
torsion, then $\mathcal{E}$ specializes to a (unique) closed point in
$\mathscr{U}$, say $\mathcal{E}_0 \in |\mathscr{U}|$, which still has torsion as
torsion-freeness is an open condition in flat families. Then
$\mathrm{Aut}(\mathcal{E}_0)$ is geometrically reductive by \cite[Prop.
9.3.4]{ams}, a contradiction by \Cref{prop: torsion is non-reductive} since $r>0$.
\end{proof}
In particular, \Cref{prop: moduli spaces in coh(C)} allows us to generalize \cite[Thm. 4.1 and 2.1]{weissmann-zhang} to the stack of coherent sheaves.
\begin{corollary}
The following statements hold:
\begin{enumerate}
\item 
The open substack $\mathscr{B}un_2^{d,ss} \subset \mathscr{C}oh_2^d$ of semistable vector bundles is the unique maximal open substack that admits a separated adequate moduli space.
\item 
The open substack $\mathscr{B}un_r^{d,ss} \subset \mathscr{C}oh_r^d$ of semistable vector bundles is the unique maximal open substack that admits a schematic adequate moduli space.
\end{enumerate}
\end{corollary}
For completeness, we end this subsection by recalling the known case when $r=0$.
\begin{proposition} \label{prop: moduli of torsion sheaves}
For any positive integer $d>0$, the stack $\mathscr{C}oh_0^d$ admits a projective adequate moduli space which is isomorphic to the Hilbert scheme $\mathrm{Hilb}^d(C) \cong C^d/\!/\mathscr{S}_d$ of length $d$ subschemes of $C$.
\end{proposition}
\begin{proof}
 This is certainly well-known, and it is explained in characteristic $0$ in \cite[Ex. 4.3.6]{huybrechts-lehn}. Since we couldn't locate the precise result in the literature for arbitrary characteristic, let us sketch a quick argument. There is a Fitting support morphism $\pi: \SCoh^d_0 \to \mathrm{Hilb}^d(C)$ defined in \cite[IV, Prop. 7.8]{rydh2008families} (note that in this case the Hilbert scheme and the Chow scheme of $C$ coincide). We shall show that $\pi$ is an adequate moduli space. 

First, note that every torsion coherent sheaf on $C$ is globally generated and has trivial first cohomology, which implies easily that $\SCoh_0^d$ is of finite type in arbitrary characteristic (cf. \cite[Cor. 1.7.7]{huybrechts-lehn}). The GIT construction of the moduli of Gieseker semistable sheaves of pure dimension $0$ in \cite[\S 1]{simpson_moduli_repn_1} applies to show that $\SCoh_0^d$ admits a projective adequate moduli space $\SCoh_0^d \to M_0^d$. Consider the morphism $f: M_0^d \to \mathrm{Hilb}^d(C)$ induced by $\pi$. The restriction of $\pi: \SCoh^d_0 \to \mathrm{Hilb}^d(C)$ over the open subscheme $\mathrm{Hilb}^d(C)_{\text{\'et}} \subset \mathrm{Hilb}^d(C)$ parametrizing distinct points in $C$ is a $\mathbb{G}_m^d$-gerbe, and hence it follows that the morphism $f: M_0^d \to \mathrm{Hilb}^d(C)$ is an isomorphism over the open subscheme $\mathrm{Hilb}^d(C)_{\text{\'et}} \subset \mathrm{Hilb}^d(C)$. 

Every torsion sheaf $\mathcal{E} \in \SCoh^d_0(k)$ is S-equivalent to the structure sheaf $\mathcal{O}_Z$ of a unique closed subscheme $Z \subset C$. Indeed, $\mathcal{E}$ is S-equivalent to a unique sheaf of the form $\mathcal{E}_0 = \oplus_{i=1}^\ell (\cO_{C}/\mathfrak{m}_{x_i})^{\oplus n_i}$ where $x_1,\ldots,x_\ell \in C(k)$ are closed points with corresponding ideal sheaves $\mathfrak{m}_{x_i}$. On the other hand, $\mathcal{E}_0$ is S-equivalent to a unique structure sheaf of a closed subscheme $Z \subset C$, namely, the one cut out by the ideal sheaf $\prod_{i=1}^\ell \mathfrak{m}_{x_i}^{n_i} \subset \cO_C$. We conclude that the morphism $f:M_0^d \to \mathrm{Hilb}^d(C)$ is bijective on $k$-points, and hence it is finite.

By \cite[Cor. A.5]{Hoffmann-Ext-stack}, the stack $\SCoh^d_0$ is integral, and hence $M^d_0$ is an integral scheme \cite[Prop. 5.4.1]{ams}. We conclude that $f: M^d_0 \to \mathrm{Hilb}^d(C)$ is a finite birational morphism between integral schemes. Since $\mathrm{Hilb}^d(C)$ is smooth, it follows that $f: M^d_0 \to \mathrm{Hilb}^d(C)$ is an isomorphism.
\end{proof}

\subsection{Reduction to semistable bundles}
\label{section:reduction to semistable vector bundles}
In this subsection we show that any substack $\mathscr{U} \subset \Bun_r^d$ with a universally closed adequate moduli space is contained in the semistable locus $\Bun_r^{d,ss}$. We shall argue by contradiction by using the following consequence of
Lange's conjecture.

\begin{lemma}
\label{lemma:general-bundles-have-homs-coprime-case}
    Let $C$ be a smooth projective connected curve of genus at least $2$ over $k$.
    Let $(d_1,r_1)$ and $(d_2,r_2)$ be pairs of coprime integers, $r_i > 0$,
    such that $d_1/r_1 < d_2/r_2$.
    Then there exist stable vector bundles $\cV_i$ of rank $r_i$ and degree $d_i$, $i=1,2$,
    such that $\Hom(\cV_1,\cV_2)\neq 0$.
\end{lemma}

\begin{proof}
    We first show the lemma in characteristic $0$. We make use of Lange's conjecture (proven in characteristic $0$ in \cite{langesconjecture-1}):
    for general bundles $\cE_1,\cE_2$ such
    that $\mu(\cE_1)<\mu(\cE_2)$ the general extension 
    $0\to \cE_1\to \cE\to \cE_2\to 0$ is stable. 
    
    If $r_1<r_2$, then we can obtain such a stable extension
    for a stable vector bundle $\cE_1$ of rank $r_1$ and degree $d_1$ and $\cE$
    of rank $r_2$ and degree $d_2$
    (since we can make sense of the correct rank $r_2 - r_1$ and degree
    $d_2-d_1$ of $\cE_2$).
    Then $\cV_1:=\cE_1$  and $\cV_2:=\cE$ is the desired example of stable 
    bundles of rank
    $r_1$ (resp. $r_2$) and degree $d_1$ (resp. $d_2$) with a non-trivial morphism
    $\cV_1\to \cV_2$.

    The case $r_1>r_2$ is similar to the case $r_1 < r_2$,
    where the extension obtained from Lange's conjecture 
    takes the role of $\cV_1$ and the quotient the role of $\cV_2$.

    If $r_1=r_2 =: r$, then $d_1<d_2$. If $r=1$, then the statement follows immediately via considering the ideal cutting out a point.
    If $r\geq 2$, then we can assume via twisting both bundles
    by a line bundle that $0\leq d_2<r$.
    From the coprime assumption on the pairs $(d_i,r_i)$,
    we also find $d_1 \neq 0$ and $d_2\neq 0$. 
    Under these assumption a computation verifies $d_1/r < (d_2 - 1)/(r-1) < d_2/r$.
    Then we can apply the above discussion first to
    obtain a surjection of stable vector bundles $\cV_1\twoheadrightarrow \cW$
    and then to obtain a monomorphism of stable vector bundles $\cW'\hookrightarrow \cV_2$,
    where $\cV_1$ is of rank $r$ and degree $d_1$, 
    $\cW$ and $\cW'$ are of rank $r-1$ and degree $d_2 -1$,
    and $\cV_2$ is of rank $r$ and degree $d_2$.
    Since Lange's conjecture holds for the general bundles,
    we can vary $\cW$ and $\cW'$
    (which also may change $\cV_1$ and $\cV_2$) 
    to be the same and find the desired non-zero morphism
    as the composition $\cV_1\twoheadrightarrow \cW \hookrightarrow \cV_2$.
    This concludes the remaining case in characteristic $0$.

    In positive characteristic, since the deformation theory of smooth 
    curves is unobstructed
    \cite[\href{https://stacks.math.columbia.edu/tag/0E84}{Tag 0E84}]{sp},
    we may choose a family
    $\mathcal{C}\to \Spec(R)$ of smooth projective
    curves over a complete discrete valuation ring $R$ of mixed characteristic
    such that the special fiber is isomorphic to $C$.
    Consider the relative moduli space $M^{d_i,s}_{r_i} \to \Spec(R)$
    of stable vector bundles of rank $r_i$ and degree $d_i$ for $i=1,2$ on the family of curves $\mathcal{C}/\Spec(R)$.
    By the coprimality assumption on $(d_i,r_i)$ we have that $M^{d_i,s}_{r_i}$
    is proper over $\Spec(R)$. Consider the locus $Z$ 
    in $M^{d_1,s}_{r_1}\times M^{d_2,s}_{r_2}$ given by pairs
    of stable vector bundles $(\cV_1,\cV_2)$ with $\Hom(\cV_1,\cV_2)\neq 0$.
    By upper-semicontinuity of dimensions of global sections, $Z$ is closed
    in $M^{d_1,s}_{r_1}\times M^{d_2,s}_{r_2}$.
    In particular, $Z$ is proper over $\Spec(R)$.
    
    Let $\eta\to \Spec(R)$ be a geometric point over the generic point.
    Then $Z_{\eta}$ is the locus of pairs of stable vector bundles $(\cV_1,\cV_2)$
    of rank $r_i$ and degree $d_i$ for $i=1,2$,
    on $\mathcal{C}_{\eta}$ admitting a non-trivial morphism between them.
    By the above discussion $Z_{\eta}$ is non-empty, which implies that $Z$ is
    non-empty as well. Thus, $Z\to \Spec(R)$ is surjective and
    we win by considering the special fiber.
\end{proof}

\begin{lemma}
    \label{lemma:general-bundles-have-homs}
    Let $C$ be a smooth projective connected curve of genus at least $2$ over $k$.
    Let $d_1,r_1, d_2, r_2$ be integers such that $r_1,r_2\geq 1$ and 
    $\mu_1:=d_1/r_1 \leq d_2/r_2=:\mu_2$.
    Then there exist semistable vector bundles $\cV_i$
    of rank $r_i$ and degree $d_i$, $i=1,2$,
    such that for any semistable vector bundle $\cW_1$ which is S-equivalent to $\cV_1$ and any semistable vector bundle $\cW_2$ which is S-equivalent to  $\cV_2$,
    we have $\Hom(\cW_1,\cW_2)\neq 0$.
\end{lemma}

\begin{proof}
    Consider $d'_i:=d_i/\mathrm{gcd}(d_i,r_i)$ and $r'_i:=r_i/\mathrm{gcd}(d_i,r_i), i=1,2$.
    If $\mu_1=\mu_2$, then we obtain the desired bundles
    $\cV_1$ and $\cV_2$ via direct sums of the same stable vector bundle
    of rank $r_1'$ and degree $d_1'$.
    If $\mu_1 < \mu_2$, then there exist stable vector bundles
    $\cV_1'$ and $\cV_2'$ of rank $r'_1$ and degree $d'_1$ 
    (resp. $r'_2$ and $d'_2$) together with a non-zero morphism
    $\cV_1'\to \cV_2'$ by \Cref{lemma:general-bundles-have-homs-coprime-case}.
    Setting $\cV_i:=(\cV'_i)^{\oplus \mathrm{gcd}(d_i,r_i)}$ we obtain the desired bundles.
    Indeed, every bundle $\cW_1$ which is S-equivalent to $\cV_1$ has at least one copy of $\cV_1'$ as a quotient
    and similarly every bundle $\cW_2$ which is S-equivalent to $\cV_2$ has at least one copy
    of $\cV_2'$ as a subbundle.
\end{proof}

\begin{proposition} \label{prop: open substack with universally closed ams is contained in semistable locus}
    Let $C$ be a smooth projective connected curve of genus $g_C \geq 2$ over $k$. Fix a choice of rank $r$ and degree $d$. Let $\mathscr{U} \subset \Bun_{r}^{d}$ be an open substack that admits a universally closed adequate moduli space. Then, $\mathscr{U} \subset \Bun_r^{d,ss}$.
\end{proposition}
\begin{proof}
Suppose for the sake of contradiction that $\mathscr{U}$ contains an unstable vector bundle. Note that the moduli space $U$ of $\mathscr{U}$ is universally closed over $k$. By \cite[\href{https://stacks.math.columbia.edu/tag/04XW}{Tag 04XW}]{sp}, it follows that $U$ is quasi-compact, and hence $\mathscr{U}$ is also quasi-compact. Consider the Harder-Narsimhan stratification of $\Bun_r^d$ by locally closed substacks (see for example \cite[Thm. 4.1]{nitsure_hn}). The quasi-compactness of $\mathscr{U}$ implies that there are finitely many strata that intersect $\mathscr{U}$. Let $\mathscr{Y}$ denote a Harder-Narasimhan stratum of $\Bun_r^d$ such that $\mathscr{U} \cap \mathscr{Y}$ is closed (note that $\mathscr{Y}$ is not the open stratum of semistable vector bundles, since $\mathscr{U}$ contains unstable vector bundles). The intersection $\mathscr{U} \cap \mathscr{Y}$ is closed in $\mathscr{U}$, so it is universally closed. On the other hand, we have that $\mathscr{U} \cap \mathscr{Y}$ is open in $\mathscr{Y}$. 
    
    The Harder-Narasimhan stratum $\mathscr{Y}$ is associated to a tuple of
    pairs of integers $(d_1, r_1), (d_2, r_2), \ldots, (d_n, r_n)$ such that
    $r_i>0$ for all $i$ and we have the strict inequalities $d_1/r_1<d_2/r_2<\ldots<d_n/r_n$. The stratum $\mathscr{Y}$ is
    the stack that parametrizes vector bundles $\cE$ on $C$ equipped with a
    filtration by subbundles $0 = \cE_{n+1} \subset \cE_{n} \subset \ldots
    \subset \cE_{2} \subset \cE_1 = \cE$ such that $\cE_i/\cE_{i+1} \in
    \Bun_{r_i}^{d_i,ss}$. In particular we have that $\mathscr{Y}$ is an open
    substack of the stack $\Bun_P$ of $P$-bundles, where $P \subset \GL_r$ is a
    parabolic subgroup. Furthermore, it follows from the description of the
    connected components of $\Bun_P$ (cf. \cite[\S 5]{hoffmann-moduli-stacks}) that
    $\mathscr{Y}$ is contained as an open inside a connected (irreducible)
    component $\Bun_P^{\vec{d}} \subset \Bun_P$ of the smooth stack $\Bun_P$.
    Hence, $\mathscr{Y}$ is an integral smooth stack. 
    
    There is a surjective morphism $\mathscr{Y} \to \prod_{i=1}^n \Bun_{r_i}^{d_i,ss}$
    that sends a filtered bundle $\cE_{\bullet}$ 
    to the tuple $(\cE_i/\cE_{i+1})_{i=1}^n$ of components of the associated graded bundle.
    Consider the surjective composition 
    $g: \mathscr{Y} \to \prod_{i=1}^n \Bun_{r_i}^{d_i,ss} \to \prod_{i=1}^n M_{r_i}^{d_i}$
    to a product of moduli spaces of bundles. 
    This is a dominant morphism from the smooth integral stack $\mathscr{Y}$. 
    In particular, the restriction to the open substack 
    $f: \mathscr{U} \cap \mathscr{Y} \to \prod_{i=1}^n M_{r_i}^{d_i}$ remains dominant.
    Since $\mathscr{U} \cap \mathscr{Y}$ is universally closed and 
    $\prod_{i=1}^n M_{r_i}^{d_i}$ is separated, 
    it follows that $f$ is universally closed.
    Hence, $f: \mathscr{U} \cap \mathscr{Y} \to \prod_{i=1}^n M_{r_i}^{d_i}$ is
    surjective.

    By \Cref{lemma:general-bundles-have-homs} above, there exist semistable
    bundles $\cV_1 \in \Bun_{r_1}^{d_1,ss}(k)$ and $\cV_n \in
    \Bun_{r_n}^{d_n,ss}(k)$ such that $\Hom(\mathcal{W}_1, \mathcal{W}_n) \neq 0$ for all
    bundles $\mathcal{W}_1,\mathcal{W}_n$ such that $\mathcal{W}_1$ 
    is S-equivalent to $\mathcal{W}_1$ and $\mathcal{W}_n$ is S-equivalent to 
    $\mathcal{W}_n$. 
    Choose arbitrary vector bundles $\cV_i \in \Bun_{r_i}^{d_i}(k)$ for $i \neq
    1,n$ and consider the image $x$ of $(\cV_1, \ldots, \cV_n)$ in the product
    of moduli spaces $\prod_{i=1}^n M_{r_i}^{d_i}$. 
    By construction, the preimage $f^{-1}(x) \subset \mathscr{U}$ is a nonempty
    closed substack of $\mathscr{U}$.

    Choose a closed $k$-point $\cE$ of $f^{-1}(x)$. Since $x$ is also a closed
    point in $\mathscr{U}$ and $\mathscr{U}$ admits an adequate moduli space, it
    follows that $\Aut(\cE)$ is geometrically reductive \cite[Prop. 9.2.4]{ams}.
    Since $\cE  \in \pi^{-1}(x)$, the Harder-Narasimhan filtration of $\cE$ is
    of the form $0 = \cE_{n+1} \subset \cE_{n} \subset \ldots \subset \cE_{2}
    \subset \cE_1 = \cE$, where $\cE_n$ is S-equivalent to $\cV_n$ and
    $\cE/\cE_2$ is S-equivalent to $\cV_1$. 
    
    There is a smooth unipotent subgroup
    $H := \mathrm{Id}_{\mathcal{E}} + \Hom(\cE/\cE_2,
    \cE_n) \subset \Aut(\cE)$,
    where we view a morphism $\cE/\cE_2 \to \cE_n$
    as an nilpotent endomorphism of $\cE$ by considering the composition 
    $\cE \twoheadrightarrow \cE/\cE_2 \to \cE_n \hookrightarrow \cE$.
    By construction we know that
    $\Hom(\cE/\cE_2, \cE_n) \neq 0$,
    and hence $H$ has positive dimension.  
    Furthermore, since the automorphisms of $\cE$ preserve
    the Harder-Narasimhan filtration, 
    it follows that $H$ is a normal subgroup of $\Aut(\cE)$.
    Therefore $\Aut(\cE)$ is not geometrically reductive, a contradiction.
\end{proof}

\begin{corollary} \label{coroll: simple bundles not universally closed}
    Let $C$ be a smooth projective connected curve of genus at least $4$ over $k$. Fix a rank $r \geq 2$ and degree $d$. Then, the open substack $\Bun_r^{d,\mathrm{simp}} \subset \Bun_r^d$ of simple vector bundles is not universally closed.
\end{corollary}
\begin{proof}
    Note that $\Bun_r^{d,\mathrm{simp}}$ admits a good moduli space which is obtained by rigidifying the central $\mathbb{G}_m$-automorphisms. By the assumption that the genus is at least $4$ and $r \geq 2$, there are simple unstable vector bundles \cite[Rmk. 2.17]{weissmann-zhang}. Therefore, it follows from \Cref{prop: open substack with universally closed ams is contained in semistable locus} that $\Bun_r^{d, \mathrm{simp}}$ is not universally closed.
\end{proof}

\begin{remark}
    We note that if $\gcd(d,r)=1$, then the good moduli space $M$ of the quasi-compact
    stack $\Bun_r^{d,\mathrm{simp}}$ has an open dense projective subspace
    $M_r^d \subset M$. However, $M$ is not universally closed.
\end{remark}

\subsection{Luna's stratification for the moduli of bundles}
\label{section: lunas stratitifcation for the moduli of bundles}
Recall that we denote by $\pi: \Bun_r^{d,ss} \to M_r^d$ the adequate moduli
space of semistable vector bundles of rank $r$ and degree $d$. In this
subsection, we introduce the Luna stratification for $M_r^d$. Although this
stratification can be defined for smooth quotient stacks admitting a good moduli space
\citep{luna-slices-etale, drezet-luna-slice, brion_luna_survey}, the moduli of
semistable vector bundles permits a more explicit and direct construction. 
\begin{definition}[Luna strata]
For any tuple $\lambda=((r_i,n_i,d_i))_{i=1}^\ell$ of triples of integers with $r_i, n_i > 0$ satisfying 
\[
\sum_{i=1}^\ell r_in_i=r \text{ and } rd_i=r_id \text{ for each } i,
\]
we consider the associated direct sum map
\[
\oplus_{\lambda}: \prod_{i=1}^\ell M_{r_i}^{d_i} \to M_r^d \text{ mapping } (\mathcal{E}_i)_{i=1}^\ell \mapsto \oplus_{i=1}^\ell \mathcal{E}_i^{\oplus n_i}.
\]
Let $\overline{S}_\lambda:=\mathrm{Im}(\oplus_{\lambda}) \subset M_r^d$ be the closed image of $\oplus_{\lambda}$. The \emph{Luna stratum} $S_\lambda$ of $M_r^d$ associated to the tuple $\lambda$ is
\[
S_\lambda:=\overline{S}_\lambda-\bigcup_{\overline{S}_\mu \subsetneq \overline{S}_\lambda} \overline{S}_\mu.
\]
\end{definition}

\begin{remark}
    Under our convention, we have that several tuples $\lambda$ can yield the same Luna stratum. If one wants to have a bijective indexing set, then we should impose that the pairs $(r_i, n_i)_{i=1}^{\ell}$ are listed in the weakly increasing lexicographic order. For simplicity we don't worry about this.
\end{remark}

\begin{example}
    If $\gcd(d,r)=1$, then there is only one tuple $\lambda$ satisfying the conditions above, i.e., $\ell=1$ and $\lambda=(r,1,d)$. In this case $S_\lambda=M_r^{d}$.
\end{example}

\begin{definition}[Minimal Luna stratum]
    For given rank $r$ and degree $d$, let 
    $r_0:=r/\gcd(d,r)$ and $d_0:=d/\gcd(d,r)$. 
    The Luna stratum of $M_r^d$ associated to the tuple $\lambda_{\min} :=(r_0,\gcd(d,r),d_0)$ with $\ell=1$
    is called the \emph{minimal} Luna stratum of $M_r^d$, denoted by $S_{\min}$.
    In this case $S_{\min}$ is the closed image of the morphism
    $\oplus_{\lambda_{\min}}: M_{r_0}^{d_0} \to M_r^d$ given by
    $\mathcal{E}_0 \mapsto \mathcal{E}_0^{\oplus \gcd(d,r)}$.
\end{definition} 

\begin{example}
In rank $3$ and degree $0$, there are five Luna strata corresponding to the following tuples $\lambda$ and associated direct sum maps:

\begin{table}[htb]
\centering
\begin{tabular}{c|c|c|c|c}
 & $\ell$ & $\lambda=((r_i,n_i,d_i))$ & $\oplus_{\lambda}$ & mapping \\
\hline
$\lambda_1$ & 1 & $(1,3,0)$ & $M_1^0 \to M_3^0$ & $\cL \mapsto \cL \oplus \cL \oplus \cL$ \\
$\lambda_2$ & 2 & $(1,1,0);(1,2,0)$ & $M_1^0 \times M_1^0 \to M_3^0$ & $(\cL_1,\cL_2) \mapsto \cL_1 \oplus \cL_2 \oplus \cL_2$ \\
$\lambda_3$ & 3 & $(1,1,0)^3$ & $M_1^0 \times M_1^0 \times M_1^0 \to M_3^0$ & $(\cL_1,\cL_2,\cL_3) \mapsto \cL_1 \oplus \cL_2 \oplus \cL_3$ \\
$\lambda_4$ & 2 & $(1,1,0);(2,1,0)$ & $M_1^0 \times M_2^0 \to M_3^0$ & $(\cL,\cF) \mapsto \cL \oplus \cF$ \\
$\lambda_5$ & 1 & $(3,1,0)$ & $M_3^0 \to M_3^0$ & $\cE \mapsto \cE$
\end{tabular}
\end{table}

Then $S_{\lambda_1}$ is the minimal stratum of $M_3^0$ and $\overline{S}_{\lambda_1} \subsetneq \overline{S}_{\lambda_2} \subsetneq \overline{S}_{\lambda_3} \subsetneq \overline{S}_{\lambda_4} \subsetneq \overline{S}_{\lambda_5}$. All Luna strata of $M_3^0$ are given by
\begin{table}[htb]
\centering
\begin{tabular}{c|c|c|c}
 & $\overline{S}_\lambda$ & $S_\lambda$ & Stabilizer of closed points in $S_\lambda$ \\
\hline
$\lambda_1$ & $\{\cL \oplus \cL \oplus \cL\}$ & $\{\cL \oplus \cL \oplus \cL\}$ & $\mathrm{GL}_3$  \\
$\lambda_2$ & $\{\cL_1 \oplus \cL_2 \oplus \cL_2\}$ & $\cL_1 \ncong \cL_2$ & $\mathbb{G}_m \times \mathrm{GL}_2$  \\
$\lambda_3$ & $\{\cL_1 \oplus \cL_2 \oplus \cL_3\}$ & $\cL_i \ncong \cL_j$ if $i \neq j$ & $\mathbb{G}_m^3$  \\
$\lambda_4$ & $\{\cL \oplus \cF\}$ & $\cF \in M_2^{0,s}$ & $\mathbb{G}_m^2$ \\
$\lambda_5$ & $\{\cE\}$ & $\cE \in M_3^{0,s}$ & $\mathbb{G}_m$ 
\end{tabular}
\end{table}

\end{example}

\begin{lemma} \label{lemma: topological properties of Luna strata}
    Let $\lambda = (r_i, n_i, d_i)_{i=1}^{\ell}$ be a tuple and let $S_\lambda \hookrightarrow M_r^d$ be the corresponding Luna stratum. Then, the following hold:
    \begin{enumerate}[(1)]
        \item $S_\lambda$ is irreducible and locally closed in $M_r^d$, 
        and $\overline{S}_{\lambda}$ is the closure of $S_{\lambda}$ in $M^d_r$.
        \item the closure $\overline{S}_\lambda$ contains the minimal stratum $S_{\min}$.
        \item for any geometric point $x \in |\pi^{-1}(S_{\lambda})| \subset |\Bun_r^{d,ss}|$ which is closed in its corresponding $\pi$-fiber, we have $\Aut(x) \cong\prod_{i=1}^\ell \mathrm{GL}_{n_i,k(x)}$.
    \end{enumerate}
\end{lemma}
\begin{proof}
As a morphism of projective varieties, the image $\overline{S}_\lambda \subset
M_r^d$ of $\oplus_\lambda$ is irreducible and closed. It follows then from
definition that $S_\lambda \subset \overline{S}_\lambda$ is open and dense, 
and hence it is irreducible and locally closed in $M_r^d$. 
This proves part (1). 
The morphism $\oplus_{\lambda_{\min}}$ factors through $\oplus_\lambda$
\[
\begin{tikzcd}
M_{r_0}^{d_0} \ar[r,"{\oplus_{\lambda_{\min}}}"] \ar[d,"h"'] & M_r^d, \\
\prod_{i=1}^\ell M_{r_i}^{d_i} \ar[ur,"{\oplus_\lambda}"']
\end{tikzcd}
\]
where $h(\mathcal{E}_0)=(\mathcal{E}_0^{\oplus r_i\gcd(d,r)/r})_{i=1}^\ell$, which shows part (2).

For part (3), note that by the definition of $S_{\lambda}$ we have that $x$ corresponds to a polystable vector bundle $\mathcal{E}$ on $C_{k(x)}$ of the form $\mathcal{E} = \oplus_{i=1}^{\ell} \mathcal{E}_i^{\oplus n_i}$ for some stable vector bundles $\mathcal{E}_i \in \Bun^{d_i,s}_{r_i}(k(x))$. It follows that $\mathrm{Aut}(\mathcal{E}) \cong \prod_{i=1}^\ell
\mathrm{GL}_{n_i,k(x)}$.
\end{proof}
We can equip each Luna stratum $S_\lambda$ with the reduced subscheme structure, so it becomes a locally closed subvariety of $M_r^d$. We have the Luna stratification
\[
M_r^d=\bigsqcup_\lambda S_\lambda.
\]

\begin{proposition} \label{prop: nonclosed points over minimal stratum don't have reductive autos}
    Let $x$ be a $k$-point in the minimal Luna stratum $S_{\min} \subset M_r^d$. Let $\mathcal{E}$ be a $k$-point in the fiber $\pi^{-1}(x) \subset \Bun_r^{d,ss}$. Then $\Aut(\mathcal{E})$ is geometrically reductive if and only if $\mathcal{E}$ is the unique closed point in $\pi^{-1}(x)$.
\end{proposition}
\begin{proof}
The if part is clear. Suppose that $\Aut(\mathcal{E})$ is geometrically reductive. The unique closed point in $\pi^{-1}(x)$ is of the form $\mathcal{E}_0^{\oplus \gcd(d,r)}$ for some stable vector bundle $\mathcal{E}_0$ of rank $r_0$ and degree $d_0$. Since $\mathcal{E}$ is S-equivalent to $\mathcal{E}_0^{\oplus \gcd(d,r)}$, we conclude by \Cref{lemma: geometrically reductive automorphism S-equivalent to direct sum of same stable bundle}.
\end{proof}

\subsection{Substacks of the semistable locus}
\label{section:substacks of the semistable locus}

In this subsection we complete the proof of 
\Cref{thm: proper moduli space is semistable locus}. 
We begin with a technical lemma regarding the interaction of
substacks of $\Bun_r^{d,ss}$ with the Luna strata. 
We keep the notation $\pi: \Bun_r^{d,ss} \to M_r^d$
for the adequate moduli space.
\begin{lemma} \label{lemma: luna strata and open substacks of bunss}
    Let $C$ be a smooth projective connected curve of genus at least $2$ over $k$. 
    Fix a choice of rank $r$ and degree $d$. 
    Let $\mathscr{U} \subset \Bun_{r}^{d,ss}$ 
    be an open substack which is universally closed. 
    Let $\mathscr{Z}= \Bun_{r}^{d,ss} \setminus \mathscr{U}$
    denote the reduced closed complement. 
    Then, $\pi(|\mathscr{Z}|) \subset |M_r^d|$ is a union of
    Luna strata in $|M_r^d|$.
\end{lemma}
\begin{proof}
    Let $S \hookrightarrow M_r^d$ be a Luna stratum. Set $\mathscr{Z}_S = \mathscr{Z} \cap \pi^{-1}(S)$. Suppose that $\pi(|\mathscr{Z}_S|) \neq |S|$. We need to show that $\pi(|\mathscr{Z}_S|)$ is empty. For the sake of contradiction, suppose that this is not the case. Then there exists a $k$-point $x \in S(k)$ which lies in the image $\pi(|\mathscr{Z}_S|)$. Note that the fiber $\pi^{-1}(x) \subset \Bun_{r}^{d,ss}$ has a unique closed $k$-point $\widetilde{x}$, and every point of $\pi^{-1}(x)$ specializes to $\widetilde{x}$. Since $\mathscr{Z}_S \cap \pi^{-1}(x) \subset \pi^{-1}(x)$ is closed, to conclude the contradiction it suffices to show that $\widetilde{x} \notin \mathscr{Z}_S(k)$ (then it would follow that $\mathscr{Z}_S \cap \pi^{-1}(x)$ is empty, contradicting $x$ being in the image of $\mathscr{Z}_S$). If we set $\mathscr{W} = \pi^{-1}(S) \cap \mathscr{U}$, then $\mathscr{W} \subset \pi^{-1}(S)$ is an open substack with closed complement $\mathscr{Z}_S$. We need to show that $\widetilde{x} \in \mathscr{W}(k)$.
    
    Note that $\pi^{-1}(S) \to S$ is universally closed, and hence the restriction to the closed substack $\pi: \mathscr{Z}_S \to S$ is also universally closed. We conclude that $\pi(|\mathscr{Z}_S|) \subset |S|$ is closed. Since $S$ is irreducible by \Cref{lemma: topological properties of Luna strata}, our assumption that $\pi(|\mathscr{Z}_S|) \neq |S|$ implies that there is an open dense subscheme $V \subset S$ such that $|V| \cap \pi(|\mathscr{Z}_S|)$ is empty. Choose a complete discrete valuation ring $R$ and a morphism $\Spec(R) \to S$ such that the generic point $\eta \in \Spec(R)$ maps to a point $v$ of $V \subset S$ and the special point $s \in \Spec(R)$ maps to $x$. 
    
     By the construction of $V$, we have $\pi^{-1}(V) \cap \mathscr{Z}_S =
     \emptyset$, and hence the unique closed point $\widetilde{v}$ of the fiber
     $\pi^{-1}(v)$ belongs to $\mathscr{W}$. Up to replacing $R$ with a finite
     discrete valuation ring extension, we may suppose that the inclusion of the
     generic point $\eta \to \Spec(R)$ lifts to a section $\eta \to \pi^{-1}(v)
     \subset \mathscr{W}$ with image the closed point $\widetilde{v}$ in
     $\pi^{-1}(v)$. Since $\mathscr{U}$ is universally closed over $k$ and
     $M_r^d$ is separated, it follows that $\pi: \mathscr{U} \to M_r^d$ is
     universally closed. Therefore $\pi: \mathscr{W} = \pi^{-1}(S) \cap
     \mathscr{U} \to S$ is universally closed. Hence, we may extend $\eta \to
     \mathscr{W}$ to a section $f:\Spec(R) \to \mathscr{W}$. 
     
     By upper-semicontinuity of dimension of automorphism groups, 
     we have that $\dim(\Aut(f(s)) \geq \dim(\Aut(f(\eta)) = \dim(\Aut(\widetilde{v}))$.
     Since $\widetilde{v}$ and $\widetilde{x}$ are points lying over the stratum
     $S$ which are closed in their $\pi$-fibers, we have
     $\dim(\Aut(\widetilde{v})) = \dim(\Aut(\widetilde{x}))$, see \Cref{lemma: topological properties of Luna strata} (3). It follows that
     $\dim(\Aut(f(s))) \geq \dim(\Aut(\widetilde{v})) =
     \dim(\Aut(\widetilde{x}))$. By \cite[Prop. 9.1]{alper-good-moduli} (this is stated for good moduli spaces but the arguments work equally for adequate moduli spaces), we
     conclude that $f(s) = \widetilde{x}$ as both map to the same point
     in $S$. 
     We have therefore shown that $\widetilde{x} \in \mathscr{W}$, as desired.
\end{proof}

\begin{remark}
    A similar argument using $\Theta$-completeness shows that in \Cref{lemma: luna strata and open substacks of bunss} we can replace the hypothesis that $\mathscr{U}$ is universally closed with that $\mathscr{U}$ has an adequate moduli space. We will not make use of this statement, so we omit its proof.
\end{remark}

\begin{proposition} \label{prop: proper moduli space inside semistable locus is semistable locus}
    Let $C$ be a smooth projective connected curve of genus at least $2$ over $k$. 
    Fix a choice of rank $r$ and degree $d$. 
    Let $\mathscr{U} \subset \mathscr{B}un_r^{d,ss}$ 
    be a nonempty open substack which admits a universally closed adequate moduli space.
    Then, $\mathscr{U} = \Bun_r^{d,ss}$.
\end{proposition}
\begin{proof}
    Let us denote by $\mathscr{Z} = \Bun_r^{d,ss} \setminus \mathscr{U}$ the
    reduced closed complement of $\mathscr{U}$. 
    We need to show that $\mathscr{Z}$ is empty, 
    or equivalently we need to prove that $\pi(|\mathscr{Z}|)$ is empty. 
    Assume for the sake of contradiction that $\pi(|\mathscr{Z}|)$ is nonempty. 
    By \Cref{lemma: luna strata and open substacks of bunss}, 
    we know that the closed subset $\pi(|\mathscr{Z}|) \subset M_r^d$
    is a union of Luna strata. 
    Since the closure of every Luna stratum contains 
    the minimal stratum $S_{\min} \subset M_r^d$, 
    see \Cref{lemma: topological properties of Luna strata}, 
    it follows that $S_{\min} \subset \pi(|\mathscr{Z}|)$. 
    We will derive a contradiction by showing that this is not the case.
    
    Since the morphism $\mathscr{U} \to M^d_r$ is dominant and is universally
    closed by assumption, it follows that $\mathscr{U} \to M_r^d$ is
    surjective. Choose a closed $k$-point $\mathcal{E}$  of $\mathscr{U}$ 
    which lies in the closed preimage $\pi^{-1}(S_{\min})$ 
    of the minimal Luna stratum $S_{\min} \subset M_r^d$.
    Since $\mathcal{E}$ is a closed $k$-point of $\mathscr{U}$ 
    and $\mathscr{U}$ admits an adequate moduli space, 
    we know that $\Aut(\mathcal{E})$ is geometrically reductive. 
    Recall that non-closed $k$-points in $\pi^{-1}(\pi(\mathcal{E}))$ cannot have geometrically reductive automorphism groups, 
    see \Cref{prop: nonclosed points over minimal stratum don't have reductive autos},
    and hence it follows that $\mathcal{E}$ must be 
    the unique $k$-point of $\pi^{-1}(\pi(\mathcal{E}))$. 
    Since every point of $\pi^{-1}(\pi(\mathcal{E}))$ 
    is a generalization of the unique closed point $\mathcal{E}$, 
    it follows that 
    $\mathscr{U} \cap \pi^{-1}(\pi(\mathcal{E})) = \pi^{-1}(\pi(\mathcal{E}))$. 
    In other words, we have shown that the image $\pi(|\mathscr{Z}|)$ 
    does not contain the point $\pi(\mathcal{E}) \in S_{\min}(k)$, 
    thus reaching our desired contradiction.
\end{proof}

\begin{theorem} \label{thm: proper moduli space is semistable locus}
    Let $C$ be a smooth projective connected curve of genus at least $2$ over $k$. 
    Fix a choice of rank $r$ and degree $d$. 
    Suppose that $\mathscr{U} \subset \mathscr{B}un_r^d$
    is a nonempty open substack which admits a universally closed adequate moduli space. 
    Then, $\mathscr{U}$ agrees with the semistable locus $\mathscr{B}un_r^{d,ss}$.
\end{theorem}
\begin{proof}
    By \Cref{prop: open substack with universally closed ams is contained in semistable locus}, 
    we have $\mathscr{U} \subset \Bun_r^{d,ss}$. 
    Then, by \Cref{prop: proper moduli space inside semistable locus is semistable locus},
    we conclude $\mathscr{U} = \Bun_r^{d,ss}$. 
\end{proof}

\section{Exotic moduli of bundles}

In this section, we provide some results and examples that illustrate what can happen if the hypothesis that the moduli space is universally closed is removed from \Cref{thm: main thm 1 intro}. We keep the same setup as in the rest of the paper: we fix a smooth projective connected curve $C$ of genus at least $2$ over an algebraically closed field $k$, and we denote by $\Bun_r^d$ the stack of rank $r$ and degree $d$ vector bundles on $C$.

\subsection{Moduli of rank \texorpdfstring{$2$}{2} bundles}
In this subsection we show that every separated adequate moduli space of rank $2$ vector bundles is necessarily an open subscheme of the moduli of semistable rank $2$ vector bundles (see \Cref{thm: moduli spaces of rank 2 vector bundles}). In particular, all separated moduli of rank $2$ vector bundles are schemes. 

For the proof of the main theorem of this subsection, we employ the Luna stratification for the moduli space $M_2^{d}$, where $d$ is even. 
In this case the stratification consists of three pieces 
$M_2^{d}=M_2^{d,s} \sqcup M_2^{d,=} \sqcup M_2^{d,\neq}$, where
\begin{gather*}
M_2^{d,=}=\{[\mathcal{L} \oplus \mathcal{L}]: \mathcal{L} \in \pic^{d/2}(C)\} \cong \pic^{d/2}(C), \\
M_2^{d,\neq}=\{[\mathcal{L}_1 \oplus \mathcal{L}_2]: 
\mathcal{L}_1 \ncong \mathcal{L}_2 \in \pic^{d/2}(C)\} \cong \left(\pic^{d/2}(C)^2 \setminus \Delta_{\pic^{d/2}(C)}\right)/\!/\mathscr{S}_2.
\end{gather*}
We shall make use of the following. 
\begin{lemma} \label{claim:irred}
   Let $d$ be an even integer. Consider the Cartesian diagram
\[
\begin{tikzcd}
\Bun_2^{d,\neq} \ar[rr,hook] \ar[d] & & \Bun_2^{d,ss} \ar[d,"\text{ams}"] \\
M_2^{d,\neq} \ar[rr,hook,"\text{loc. closed}"'] & & M_2^{d}.
\end{tikzcd}
\]
The locally closed substack $\Bun_2^{d,\neq} \hookrightarrow \Bun_2^{d,ss}$ is irreducible.
\end{lemma}
\begin{proof}
Consider the algebraic stack $\mathscr{X}:= \mathscr{E}xt((1,d/2),(1,d/2))$
parametrizing extensions
$0 \to \mathcal{L}_1 \to \mathcal{V} \to \mathcal{L}_2 \to 0$, where
$\mathcal{L}_i$ is a coherent sheaf on $C$ of rank $1$ and degree $d/2$
for $i=1,2$,
see, e.g. \cite[Appx. A]{Hoffmann-Ext-stack}. 
The $\mathscr{E}xt$-stack $\mathscr{X}$ comes equipped with two natural morphisms:
\[
\mathrm{pr}_2: \mathscr{X} \to \mathscr{C}oh_2^d, \text{ and } \mathrm{pr}_{13}: \mathscr{X} \to \mathscr{C}oh_1^{d/2} \times \mathscr{C}oh_1^{d/2}
\]
given by $\mathrm{pr}_2([\mathcal{L}_1 \to \cV \to \mathcal{L}_2]) = \cV$ and
$\mathrm{pr}_{13}([\mathcal{L}_1 \to \cV \to \mathcal{L}_2]) =
(\mathcal{L}_1,\mathcal{L}_2)$.
By definition there is a surjection 
\[
\mathrm{pr}_2: \mathrm{pr}_{13}^{-1}((\Pic^{d/2}\times \Pic^{d/2}) \setminus \Delta_{\Pic^{d/2}}) \twoheadrightarrow \Bun_2^{d,\neq},
\]
where $\Delta_{\Pic^{d/2}}$ is the (closed) image of the diagonal morphism.
This implies that $\Bun_2^{d,\neq}$ is irreducible, as 
\[
(\Pic^{d/2}\times \Pic^{d/2}) \setminus \Delta_{\Pic^{d/2}} \subset \Pic^{d/2} \times \Pic^{d/2} \subset \mathscr{C}oh_1^{d/2} \times \mathscr{C}oh_1^{d/2}
\]
is open and the $\mathscr{E}xt$-stack is irreducible by \cite[Prop. A.3 and A.4]{Hoffmann-Ext-stack}.
\end{proof}

\begin{theorem}\label{thm: moduli spaces of rank 2 vector bundles}
Let $d$ be an integer.
Suppose that $\mathscr{U} \subset \Bun^d_2$ is an open substack that admits a
separated adequate moduli space $\pi':\mathscr{U} \to U$. Then, $\mathscr{U}$ is an
open substack of $\Bun^{d,ss}_2$ that is saturated with respect to the adequate
moduli space $\pi: \Bun^{d,ss}_2 \to M_2^{d}$, and $U$ is an open subscheme of
$M_2^{d}$.
\end{theorem}
\begin{proof}
Any open substack of $\Bun_2^d$ admitting a separated adequate moduli space
is contained in the semistable locus by \cite[Thm. 4.1]{weissmann-zhang},
so we have $\mathscr{U} \subset \Bun_2^{d,ss}$. 
It suffices to show that $\mathscr{U}$ is saturated with respect
to $\pi:\Bun^{d,ss}_2\to M^d_2$ as adequate moduli spaces are preserved
by flat base change \cite[Prop. 5.2.9 (1)]{ams}. That is, we wish to show that
$\pi^{-1}(\pi(\mathscr{U}))=\mathscr{U}$. 
It suffices to show that $\mathscr{U}$ contains
the $\pi$-fiber of every $k$-point of $M^{d}_2$ which lies 
in the image of $\pi_{\mid \mathscr{U}}: \mathscr{U} \to M_2^d$.
As such a fiber contains a unique closed point, namely 
the unique polystable vector bundle in the fiber, it suffices to show that 
for a semistable vector bundle $\mathcal{V}$ lying in $\mathscr{U}(k)$
the unique polystable vector bundle S-equivalent to $\mathcal{V}$
also lies in $\mathscr{U}$.

There are three cases: 
$\mathcal{V}$ is stable, 
$\mathcal{V}$ is S-equivalent to $\mathcal{L} \oplus \mathcal{L}$ for a line bundle 
$\mathcal{L}$,  and $\mathcal{V}$ is S-equivalent to $\mathcal{L}_1\oplus \mathcal{L}_2$ 
for line bundles $\mathcal{L}_1\ncong \mathcal{L}_2$ of the same  degree. 
The first case is immediate: the $\pi$-fiber over a stable
bundle consists only of the bundle itself and $\mathscr{U}\to U$ is surjective.

If $\mathcal{V}$ is S-equivalent to $\mathcal{L} \oplus \mathcal{L}$, 
then any non-closed point $\pi^{-1}(\pi(\mathcal{V}))$ 
has non-reductive automorphism group by
\Cref{prop: nonclosed points over minimal stratum don't have reductive autos}.
However, as $\pi':\mathscr{U}\to U$ is an adequate moduli space,
every fiber of $\pi'$ contains a unique closed point with a geometrically reductive
automorphism group. Thus, we find that $\mathcal{L}\oplus\mathcal{L}$ 
is a $k$-point of $\mathscr{U}$
as the only point with geometrically reductive automorphism group in its fiber.

The remaining case when $\mathcal{V}$ is S-equivalent to
$\mathcal{L}_1\oplus \mathcal{L}_2$ with $\mathcal{L}_1\ncong \mathcal{L}_2$ 
takes more work and we need to introduce some notation. 
Note that $\mathcal{V}$ is an extension
of $\mathcal{L}_2$ by $\mathcal{L}_1$ (or vice versa), so without loss of generality
assume that we have
$0\to \mathcal{L}_1\to \mathcal{V}\to \mathcal{L}_2\to 0$. 
Also note that $\mathcal{L}_1$ and $\mathcal{L}_2$ have
degree $d':=d/2$ as $V$ is semistable; in particular $d$ is even in this case.

Similarly to \Cref{claim:irred}
we consider the Ext-stack $\mathscr{X}:=\mathscr{E}xt(1,d'),(1,d'))$ parametrizing
extensions of line bundles of degree $d'$ together
with the projections $\mathrm{pr}_{13}:\mathscr{X}\to\Pic^{d'}\times\Pic^{d'}$
and $\mathrm{pr}_2:\mathscr{X}\to \Bun^d_2$. 
Recall that $\mathrm{pr}_{13}$ is smooth \cite[Prop. A.3 (ii)]{Hoffmann-Ext-stack}.
The open preimage $\mathrm{pr}^{-1}_2(\mathscr{U})\subset \mathscr{X}$ is non-empty,
since it contains $\mathcal{V}$ considered as an extension 
of $\mathcal{L}_2$ by $\mathcal{L}_1$.
From the smoothness of $\mathrm{pr}_{13}$ we also obtain a non-empty open 
$\mathscr{W}:=\mathrm{pr}_{13}(\mathrm{pr}_2^{-1}(\mathscr{U}))$ in
$\Pic^{d'}\times \Pic^{d'}$.
Consider the automorphism $swap:\Pic^{d'}\times \Pic^{d'}\to
\Pic^{d'}\times \Pic^{d'}, (\mathcal{L},\mathcal{L}')\mapsto (\mathcal{L}',\mathcal{L})$.
As $\mathscr{W}$ is non-empty open, so is $swap^{-1}(\mathscr{W})$.
Furthermore, the product $\Pic^{d'}\times \Pic^{d'}$ is irreducible
and so the intersection $\mathscr{W}':=\mathscr{W}\cap swap^{-1}(\mathscr{W}) 
\cap (\Pic^{d'}\times \Pic^{d'}\setminus \Delta_{\Pic^{d'}})$ 
is non-empty and open, where we denote by $\Delta_{\Pic^{d'}}$ the (closed) image of the diagonal
$\Pic^{d'} \to \Pic^{d'} \times \Pic^{d'}$.
By definition $\mathscr{W}'$ 
parametrizes pairs of non-isomorphic line bundles
$(\mathcal{L},\mathcal{L}')$ of degree $d'$ such that
we have extensions of the line bundles (in both directions) lying in $\mathscr{U}$.

We claim that
the morphism $\mathscr{W}' \to \Bun_2^{d,ss}$ given by $(\mathcal{L},
\mathcal{L}') \mapsto \mathcal{L} \oplus \mathcal{L}'$ factors through
$\mathcal{U} \subset \Bun_2^{d,ss}$. Indeed, in this setting the separatedness
of the adequate moduli space $U$ implies the S-completeness of $\mathscr{U}$
\cite[Prop. 3.28]{weissmann-zhang}.
This in turn implies that $\mathscr{U}$
contains all required direct sums $\mathcal{L} \oplus \mathcal{L}'$ by 
considering the short exact sequences defining an extension
as a filtration as follows: Allowing for extensions of line bundles $\mathcal{L},\mathcal{L}'$
in both directions
yields bundles with opposite filtrations in the sense of \cite[Defn.
3.15]{weissmann-zhang}. By \cite[Prop. 3.17]{weissmann-zhang} 
we find that the associated graded $\mathcal{L}\oplus\mathcal{L}'$ of the
filtrations also lies in $\mathscr{U}$.

By \Cref{claim:irred} the stack $\Bun_2^{d,\neq}$ is irreducible. 
Denote by $\mathcal{E}_{gen} \in \Bun_2^{d,\neq}(K)$ a point defined
over an algebraically closed field $K$ that lies over the generic point, 
and let $0 \to \mathcal{V}_1 \to \mathcal{E}_{gen} \to \mathcal{V}_2 \to 0$ 
be a Jordan-H\"{o}lder filtration. 
Then $[\mathcal{V}_1 \oplus \mathcal{V}_2] \in M_2^{d,\neq}(K)$ lies over 
the generic point and hence $[\mathcal{V}_1 \oplus \mathcal{V}_2]$ must be in the image of the dominant morphism $\mathscr{W}' \hookrightarrow \Pic^{d'} \times \Pic^{d'}
\setminus \Delta_{\Pic^{d'}} \to M_2^{d,\neq}$.
By the discussion above we find that $\mathcal{V}_1\oplus\mathcal{V}_2$ also lies
in $\mathscr{U}$.

Recall that we started out with a bundle $\mathcal{V}$ which was S-equivalent 
to a direct sum of line bundles $\mathcal{L}_1\oplus \mathcal{L}_2$.
If $\mathcal{V}\cong \mathcal{L}_1 \oplus \mathcal{L}_2$ is split, 
then there is nothing to prove. 
Otherwise $0 \to \mathcal{L}_1 \to \mathcal{V} \to \mathcal{L}_2 \to 0$ is the
unique Jordan-H\"{o}lder filtration of $\mathcal{V}$. 
Since $\mathscr{U} \cap \Bun_2^{d,\neq}$ is irreducible,  there is a morphism 
$\mathrm{Spec}(R) \to \mathscr{U} \cap \Bun_2^{d,\neq} \subset \mathscr{U}$, 
where $R$ is a discrete valuation ring, 
such that $\eta \mapsto \mathcal{E}_{gen}$ and $0 \mapsto \mathcal{V}$. 
The filtration $0 \subset \mathcal{V}_1 \subset \mathcal{E}_{gen}$ 
of $\mathcal{E}_{gen}$ satisfies
$\mathcal{V}_1 \oplus \mathcal{V}_2 \in |\mathscr{U}|$ 
as we have seen above.
Thus, the morphism $\mathrm{Spec}(R) \to \mathscr{U}$ and the
filtration $0 \subset \mathcal{V}_1 \subset \mathcal{E}_{gen}$ define a morphism
$\Theta_R-\{0\} \to \mathscr{U}$ (see \cite[Lem. 3.4]{weissmann-zhang}).
This morphism extends to $\Theta_R$ since
$\mathscr{U}$ is $\Theta$-complete, see \cite[Prop. 3.28]{weissmann-zhang}. 
The image of $0$ in $\mathscr{U}$ is the
associated graded sheaf of the filtration of $\mathcal{V}$ induced from 
$0 \subset \mathcal{V}_1 \subset \mathcal{E}_{gen}$, 
which has to be the unique Jordan-H\"{o}lder filtration of $\mathcal{V}$, i.e.,
$0 \to \mathcal{L}_1 \to \mathcal{V} \to \mathcal{L}_2 \to 0.$
This shows that $\mathcal{L}_1 \oplus \mathcal{L}_2 \in |\mathscr{U}|$.

\end{proof}
\begin{example}
The conclusion of \Cref{thm: moduli spaces of rank 2 vector bundles} fails if we do not assume that the adequate moduli space $U$ is separated. For example, the open substack 
\[
\mathscr{U}:=\mathscr{B}un_2^{0,ss} \cap \mathscr{B}un_2^{0,\text{simple}} \subset \mathscr{B}un_2^{0,ss}
\]
admits a good moduli space $\mathscr{U} \to U$ via rigidifying the central $\mathbb{G}_m$-automorphisms (see \cite[Thm. 5.1.5]{acv_twisted_bundles}). We note that $U$ is non-schematic and non-separated. Indeed, choose two degree $0$ line bundles $\mathcal{L}_1 \ncong \mathcal{L}_2$ on $C$. For any non-split extensions
\[
0 \to \mathcal{L}_1 \to \mathcal{E} \to \mathcal{L}_2 \to 0 \text{ and } 0 \to \mathcal{L}_2 \to \mathcal{E}' \to \mathcal{L}_1 \to 0,
\]
it follows that $\mathcal{E},\mathcal{E}'$ are semistable and simple 
by \Cref{thm:non-split}. 
Hence they define two $k$-points in $\mathscr{U}$.
These two $k$-points are identified in $M_2^0$, 
but not identified in $U$. 
This shows that $U$ is not a scheme by 
\cite[Thm. 5.11]{distinguishing_spaces_paper}. 
Furthermore, these two $k$-points have opposite filtrations but $\mathcal{L}_1 \oplus \mathcal{L}_2 \notin \mathscr{U}(k)$, so $U$ is not separated by \cite[Prop. 3.17]{weissmann-zhang}.
\end{example}
\subsection{Exotic moduli of rank \texorpdfstring{$3$}{3} semistable bundles}
In this subsection we provide an example 
of an open substack $\mathscr{U}$ of $\Bun_3^{0,ss}$
which admits a non-schematic separated non-proper good moduli space. Let us start by defining $\mathscr{U}$.

\begin{proposition}
\label{prop: definition U}
    Let $|\mathscr{U}| \subset |\Bun_3^{0,ss}|$ denote the locus of vector
    bundles that do not admit a non-zero morphism from any degree $0$ line bundle. 
    Then $|\mathscr{U}|$ is open in $|\Bun_3^{0,ss}|$ and nonempty. 
    We denote by $\mathscr{U} \subset \Bun_3^{0,ss}$ the corresponding nonempty
    open substack of $\Bun_3^{0,ss}$.
\end{proposition}
\begin{proof}
We have $|\Bun_3^{0,s}| \subset |\mathscr{U}|$, so $|\mathscr{U}|$ is nonempty as the genus of $C$ is at least $2$. To see the openness of $|\mathscr{U}| \subset |\Bun^{0,ss}_3|$, we consider the projections:
\[
\begin{tikzcd}
\Bun^{0,ss}_3 \times C & \Bun^{0,ss}_3 \times \Pic^0(C) \times C \ar[l,"p_{13}"'] \ar[r,"p_{23}"] \ar[d,"p_{12}"] & \Pic^0(C) \times C \\
& \Bun^{0,ss}_3 \times \Pic^0(C) \ar[d,"p_1"] \\
& \Bun^{0,ss}_3
\end{tikzcd}
\]
Let $\mathcal{E}_{univ}$ (resp., $\mathcal{L}_{univ}$) be the universal family on $\Bun^{0,ss}_3 \times C$ (resp., $\Pic^0(C) \times C$). Set $\mathcal{E}:=p_{13}^*\mathcal{E}_{univ} \otimes p_{23}^*\mathcal{L}_{univ}$. By cohomology and base-change, the locus
\begin{align*}
\mathscr{Z}:&=\{s \in |\Bun^{0,ss}_3 \times \Pic^0(C)|: H^0(C_{k(s)},\mathcal{E}_s) \neq 0\} \subset |\Bun^{0,ss}_3 \times \Pic^0(C)|
\end{align*}
is closed. Note that $\Pic^0(C) \to \pic^0(C) \to \Spec(k)$ is universally closed, the projection $p_1: \Bun^{0,ss}_3 \times \Pic^0(C) \to \Bun^{0,ss}_3$ is closed. Then $|\mathscr{U}| \subset |\Bun^{0,ss}_3|$ is the open complement of the closed subset $p_1(\mathscr{Z}) \subset |\Bun^{0,ss}_3|$.
\end{proof}
\begin{lemma}\label{lemma: ss-in-U}
We keep the notation of \Cref{prop: definition U}.
Let $\mathcal{E} \in |\Bun_3^{0,ss}|$ be a geometric point. 
Then $\mathcal{E} \in |\mathscr{U}|$
if and only if it is stable or it fits into a non-split short exact sequence $0
\to \mathcal{F} \to \mathcal{E} \to \mathcal{L} \to 0$ for some stable rank $2$,
degree $0$ vector bundle $\mathcal{F}$ and degree $0$ line bundle $\mathcal{L}$.
In the latter case, 
$\mathcal{E}$ is strictly semistable and 
$0 \subset \mathcal{F} \subset \mathcal{E}$ is its unique Jordan-H\"{o}lder filtration.
\end{lemma}
\begin{proof}
Suppose $\mathcal{E} \in |\mathscr{U}|$ is strictly semistable. Since
$\mathcal{E}$ does not admit a non-zero morphism from any degree $0$ line
bundle, its Jordan-H\"{o}lder filtrations are necessarily of the non-split form $0 \to \mathcal{F} \to \mathcal{E} \to \mathcal{L} \to 0$ for some stable rank $2$, degree $0$ vector bundle $\mathcal{F}$ and degree $0$ line bundle $\mathcal{L}$.

Conversely, suppose $\mathcal{E} \in |\Bun_3^{0,ss}|$ fits into such a non-split short exact sequence $0 \to \mathcal{F} \to \mathcal{E} \to \mathcal{L} \to 0$. For any non-zero morphism $\mathcal{L}' \to \mathcal{E}$ from a degree $0$ line bundle, the composition $\mathcal{L}' \to \mathcal{E} \to \mathcal{L}$ is either an isomorphism or zero. If it is an isomorphism, then the short exact sequence splits. If it is zero, then the non-zero morphism $\mathcal{L}' \to \mathcal{E}$ factors through $\mathcal{L}' \to \mathcal{F}$, contradicting the stability of $\mathcal{F}$.
\end{proof}
\begin{corollary}\label{lemma: all-simple}
Every vector bundle in $\mathscr{U}$ is simple.
\end{corollary}
\begin{proof}
Let $\mathcal{E} \in |\mathscr{U}|$. If $\mathcal{E}$ is stable, then it is automatically simple. If $\mathcal{E}$ is strictly semistable, then by \Cref{lemma: ss-in-U} it fits into a non-split short exact sequence of the form $0 \to \mathcal{F} \to \mathcal{E} \to \mathcal{L} \to 0$ for some stable rank $2$, degree $0$ vector bundle $\mathcal{F}$ 
and degree $0$ line bundle $\mathcal{L}$. 
It follows that $\mathcal{E}$ is simple by \Cref{thm:non-split}.
\end{proof}
Therefore, we see that 
\[
\Bun_3^{0,s} \subsetneq \mathscr{U} \subsetneq \Bun_3^{2,ss} \cap \Bun_3^{0,simple}.
\]
The second inclusion is strict. Indeed, any non-split extension $0 \to
\mathcal{L} \to \mathcal{E} \to \mathcal{F} \to 0$ for some stable rank $2$,
degree $0$ vector bundle $\mathcal{F}$ and degree $0$ line bundle $\mathcal{L}$
gives a semistable, simple (by \Cref{thm:non-split}) vector bundle which is not
in $\mathscr{U}$.
\begin{proposition}\label{prop: open-in-rank-3}
The open substack $\mathscr{U} \subset \Bun_3^{0,ss}$ admits a non-schematic separated non-proper good moduli space.
\end{proposition}
\begin{proof}
By \Cref{lemma: all-simple}, all vector bundles in $\mathscr{U}$ are simple, and hence $\mathscr{U}$ admits a good moduli space $\mathscr{U} \to U$ via rigidification by the central $\mathbb{G}_m$-automorphisms \cite[Thm. 5.1.5]{acv_twisted_bundles}. The non-properness follows from \Cref{prop: proper moduli space inside semistable locus is semistable locus} (we can also directly show that the existence part of the valuative criterion for properness does not hold for $\mathscr{U}$ by mimicking the constructions in \citep[Lem 5.10 and Cor. 5.11]{weissmann-zhang}). It remains to show that $U$ is separated and that it is not a scheme.

The separatedness of $U$ is equivalent to the S-completeness of $\mathscr{U}$,
see \cite[Prop. 3.28]{weissmann-zhang}.
For open substacks of coherent sheaves, S-completeness can be checked via opposite filtrations in the following sense:
For all points $\mathcal{V} \ncong \mathcal{W} \in \mathscr{U}(K)$, where $K/k$ is a field extension, and finite $\mathbb{Z}$-graded filtrations of subsheaves
\begin{gather*}
    \mathcal{V}^{\bullet} = \left( 0 \subset \dots \subset \mathcal{V}^i \subset \mathcal{V}^{i+1} \subset \dots \subset \mathcal{V} \right), \\
    \mathcal{W}_{\bullet} = \left( \mathcal{W} \supset \dots \supset \mathcal{W}_{i+1} \supset \mathcal{W}_i \supset \dots \supset 0 \right)
\end{gather*}
such that there is an isomorphism $\mathcal{V}^{i}/\mathcal{V}^{i-1}\cong \mathcal{W}_{i+1}/\mathcal{W}_{i}$ for all $i$, then
the associated graded  sheaf  
$\oplus_{i}\mathcal{V}^{i}/\mathcal{V}^{i-1} \cong
\oplus_{i}\mathcal{W}_{i+1}/\mathcal{W}_{i}$ is a $K$-point of $\mathscr{U}$,
see \cite[Prop. 3.17]{weissmann-zhang}. Opposite filtrations of semistable
vector bundles can be refined to be their Jordan-H\"{o}lder filtrations by the
proof of \cite[Cor. 3.18]{weissmann-zhang}. Thus, the existence of opposite
filtrations of rank $3$, degree $0$ semistable vector bundles $\mathcal{V},\mathcal{W}$
requires at least one of $\mathcal{V},\mathcal{W}$ to contain a degree $0$ line
bundle. This shows that no points in $\mathscr{U}$ have opposite filtrations,
and we obtain the desired S-completeness.

Let $\Bun_3^{0,ss} \to M_3^{0}$ be the adequate moduli space and 
$f: U \to M_3^{0}$ be the induced morphism on the adequate moduli spaces. 
To show that $U$ is not a scheme, we only need to prove that $f$ 
is not an open immersion by \cite[Thm. 5.11]{distinguishing_spaces_paper}. 
To see this, let us show that $f$ is not injective on points. 
Indeed, since the genus of $C$ is at least $2$, we can choose a
stable rank $2$, degree $0$ vector bundle $\mathcal{F}$ and a degree $0$ line bundle
$\mathcal{L}$ such that $\mathrm{Ext}^1(\mathcal{L},\mathcal{F}) \neq 0$. 
Then any non-zero element in $\mathrm{Ext}^1(\mathcal{L},\mathcal{F})$ defines a
point in $\mathscr{U}$ by \Cref{lemma: ss-in-U}. 
All non-zero elements in $\mathrm{Ext}^1(\mathcal{L},\mathcal{F})$ 
are identified in $M_3^{0}$ as $[\mathcal{F} \oplus \mathcal{L}]$ 
since they are all S-equivalent, 
but they are not identified in $U$ as $\mathscr{U}$ 
contains only indecomposable vector bundles. 
This means the fiber $f^{-1}([\mathcal{F} \oplus \mathcal{L}])$ contains the
positive dimensional space $\mathrm{Ext}^1(\mathcal{L},\mathcal{F})$.

\end{proof}

\bibliographystyle{alpha}

\bibliography{bibliography}

\begin{thebibliography}{{Sta}25}

\bibitem[ACV03]{acv_twisted_bundles}
Dan Abramovich, Alessio Corti, and Angelo Vistoli.
\newblock Twisted bundles and admissible covers.
\newblock {\em Comm. Algebra}, 31(8):3547--3618, 2003.
\newblock Special issue in honor of Steven L. Kleiman.

\bibitem[AHR20]{luna-slice-for-stacks}
Jarod Alper, Jack Hall, and David Rydh.
\newblock A {L}una \'etale slice theorem for algebraic stacks.
\newblock {\em Ann. of Math. (2)}, 191(3):675--738, 2020.

\bibitem[Alp13]{alper-good-moduli}
Jarod Alper.
\newblock Good moduli spaces for {A}rtin stacks.
\newblock {\em Ann. Inst. Fourier (Grenoble)}, 63(6):2349--2402, 2013.

\bibitem[Alp14]{ams}
Jarod Alper.
\newblock Adequate moduli spaces and geometrically reductive group schemes.
\newblock {\em Algebraic Geometry}, 1(4):489--531, 2014.

\bibitem[Ati56]{MR0086358}
M.~Atiyah.
\newblock On the {K}rull-{S}chmidt theorem with application to sheaves.
\newblock {\em Bull. Soc. Math. France}, 84:307--317, 1956.

\bibitem[BKS25]{brion_luna_survey}
M.~Brion, H.~Kraft, and G.~Schwarz.
\newblock Luna's slice theorem and applications, 2025.
\newblock Available at \url{https://www-fourier.univ-grenoble-alpes.fr/~mbrion/LunaSlice.pdf}.

\bibitem[Don83]{donaldson-narasimhan-seshadri}
S.~K. Donaldson.
\newblock A new proof of a theorem of {N}arasimhan and {S}eshadri.
\newblock {\em J. Differential Geom.}, 18(2):269--277, 1983.

\bibitem[Dre04]{drezet-luna-slice}
Jean-Marc Drezet.
\newblock Luna's slice theorem and applications.
\newblock In {\em Algebraic group actions and quotients}, pages 39--89. Hindawi Publ. Corp., Cairo, 2004.

\bibitem[HL97]{huybrechts-lehn}
Daniel Huybrechts and Manfred Lehn.
\newblock {\em The geometry of moduli spaces of sheaves}, volume E31 of {\em Aspects of Mathematics}.
\newblock Friedr. Vieweg \& Sohn, Braunschweig, 1997.

\bibitem[Hof10a]{Hoffmann-Ext-stack}
Norbert Hoffmann.
\newblock Moduli stacks of vector bundles on curves and the {K}ing-{S}chofield rationality proof.
\newblock In {\em Cohomological and geometric approaches to rationality problems}, volume 282 of {\em Progr. Math.}, pages 133--148. Birkh\"{a}user Boston, Inc., Boston, MA, 2010.

\bibitem[Hof10b]{hoffmann-moduli-stacks}
Norbert Hoffmann.
\newblock On moduli stacks of {$G$}-bundles over a curve.
\newblock In {\em Affine flag manifolds and principal bundles}, Trends Math., pages 155--163. Birkh\"{a}user/Springer Basel AG, Basel, 2010.

\bibitem[HWZ24]{distinguishing_spaces_paper}
Andres~Fernandez Herrero, Dario Wei{\ss}mann, and Xucheng Zhang.
\newblock Distinguishing algebraic spaces from schemes.
\newblock \href{https://arxiv.org/pdf/2411.07169.pdf}{arXiv: 2411.07169}, 2024.

\bibitem[Lan02]{Lang-alg}
Serge Lang.
\newblock {\em Algebra}, volume 211 of {\em Graduate Texts in Mathematics}.
\newblock Springer-Verlag, New York, third edition, 2002.

\bibitem[Lun73]{luna-slices-etale}
Domingo Luna.
\newblock Slices \'etales.
\newblock In {\em Sur les groupes alg\'ebriques}, volume Tome 101 of {\em Suppl\'ement au Bull. Soc. Math. France}, pages 81--105. Soc. Math. France, Paris, 1973.

\bibitem[MFK94]{mumford_git}
D.~Mumford, J.~Fogarty, and F.~Kirwan.
\newblock {\em Geometric invariant theory}, volume~34 of {\em Ergebnisse der Mathematik und ihrer Grenzgebiete (2) [Results in Mathematics and Related Areas (2)]}.
\newblock Springer-Verlag, Berlin, third edition, 1994.

\bibitem[Mum63]{mumford-projetive-invariants}
David Mumford.
\newblock Projective invariants of projective structures and applications.
\newblock In {\em Proc. {I}nternat. {C}ongr. {M}athematicians ({S}tockholm, 1962)}, pages 526--530. Inst. Mittag-Leffler, Djursholm, 1963.

\bibitem[Nit11]{nitsure_hn}
Nitin Nitsure.
\newblock Schematic {H}arder-{N}arasimhan stratification.
\newblock {\em Internat. J. Math.}, 22(10):1365--1373, 2011.

\bibitem[NS65]{narasimhan-seshadri-unitary}
M.~S. Narasimhan and C.~S. Seshadri.
\newblock Stable and unitary vector bundles on a compact {R}iemann surface.
\newblock {\em Ann. of Math. (2)}, 82:540--567, 1965.

\bibitem[RTiB99]{langesconjecture-1}
Barbara Russo and Montserrat Teixidor~i Bigas.
\newblock On a conjecture of {L}ange.
\newblock {\em J. Algebraic Geom.}, 8(3):483--496, 1999.

\bibitem[Ryd08]{rydh2008families}
David Rydh.
\newblock {\em Families of cycles and the Chow scheme}.
\newblock KTH, 2008.
\newblock Thesis (Ph.D.)--KTH, Stockholm. pp. 218. Available at: \url{https://people.kth.se/~dary/thesis/thesis.pdf}.

\bibitem[Sch16]{indecom-Higgs}
Olivier Schiffmann.
\newblock Indecomposable vector bundles and stable {H}iggs bundles over smooth projective curves.
\newblock {\em Ann. of Math. (2)}, 183(1):297--362, 2016.

\bibitem[Ses67]{seshadri-unitary}
C.~S. Seshadri.
\newblock Space of unitary vector bundles on a compact {R}iemann surface.
\newblock {\em Ann. of Math. (2)}, 85:303--336, 1967.

\bibitem[Ses82]{seshadri}
C.~S. Seshadri.
\newblock {\em Fibr\'es vectoriels sur les courbes alg\'ebriques}, volume~96 of {\em Ast\'erisque}.
\newblock Soci\'et\'e{} Math\'ematique de France, Paris, 1982.
\newblock Notes written by J.-M. Drezet from a course at the \'Ecole Normale Sup\'erieure, June 1980.

\bibitem[Sim94]{simpson_moduli_repn_1}
Carlos~T. Simpson.
\newblock Moduli of representations of the fundamental group of a smooth projective variety. {I}.
\newblock {\em Inst. Hautes \'Etudes Sci. Publ. Math.}, 79:47--129, 1994.

\bibitem[{Sta}25]{sp}
The {Stacks Project Authors}.
\newblock \textit{Stacks Project}.
\newblock \url{https://stacks.math.columbia.edu}, 2025.

\bibitem[WZ25]{weissmann-zhang}
Dario Wei{\ss}mann and Xucheng Zhang.
\newblock A stacky approach to identifying the semistable locus of bundles.
\newblock {\em Algebr. Geom.}, 12(2):262--298, 2025.

\end{thebibliography}

\end{document}